\documentclass[12pt]{amsart}
\usepackage{amsmath,amsthm, amssymb, amsfonts, epsfig,amscd}
\usepackage[all]{xy}

\usepackage[mathscr]{eucal}
\usepackage{mathrsfs}
\usepackage{latexsym}
\usepackage{graphicx}
\newdimen\harrowlength \harrowlength=60pt
\newdimen\varrowlength \varrowlength=.618\harrowlength
\newdimen\sarrowlength \sarrowlength=\harrowlength

\newdimen\hgrid \hgrid=15pt
\newdimen\vgrid \vgrid=15pt

\newdimen\hchannel  \hchannel=0pt
\newdimen\vchannel  \vchannel=0pt
\newdimen\channelwidth \channelwidth=3pt

\theoremstyle{plain}
\newtheorem{thm}[subsection]{Theorem}
\newtheorem{lem}[subsection]{Lemma}
\newtheorem{prop}[subsection]{Proposition}
\newtheorem{cor}[subsection]{Corollary}
\newtheorem{lemma}[subsection]{Lemma}
\newtheorem{theorem}[subsection]{Theorem}

\theoremstyle{definition}
\newtheorem{rem}[subsection]{Remark}
\newtheorem{defn}[subsection]{Definition}
\newtheorem{definition}[subsection]{Definition}
\newtheorem{remark}[subsection]{Remark}

\def\L{\mathcal{L}}
\newcommand{\Z}{{\mathbb Z}}
\newcommand{\R}{{\mathbb R}}

\def\a{\alpha}
\def\t{\tau}

\def\f{\phi}
\def\l{\lambda}

\def\m{\mu}
\def\n{\nu}

\def\v{\vee}
\def\w{\wedge}

\def\Spec{\mbox{Spec}}

\newcommand{\Lij}{{{\mathcal{L}}_{ij}}}

\def\vs{\vskip}
\def\ni{\noindent}

\begin{document}
\title[Singular loci of G-H toric varieties]{Singular loci of Grassmann-Hibi toric varieties}
\address{Department of Mathematics\\ Northeastern University\\ Boston, MA 02115}
\author[J. Brown]{J. Brown}
\author[V. Lakshmibai]{V. Lakshmibai${}^{\dag}$}

\email{brown.justin1@neu.edu, lakshmibai@neu.edu}
\thanks{${}^{\dag}$ Partially supported
by NSF grant DMS-0400679 and NSA-MDA904-03-1-0034.}

\begin{abstract}
For the toric variety $X$ associated to the Bruhat poset of
Schubert varieties in the Grassmannian, we describe the singular
locus in terms of the faces of the associated polyhedral cone. We
further show that the singular locus is pure of codimension 3 in
$X$, and the generic singularities are of cone type. We also
determine the tangent cones at the maximal singularities of $X$.
These turn out to be again toric varieties. In the case of $X$
being associated to the  Bruhat poset of Schubert varieties in the
Grassmannian of $2$-planes in a $n$-dimensional vector space (over
the base field), we also prove a certain product formula relating
the multiplicities at certain singular points.
\end{abstract}
\maketitle

\section*{Introduction}
Let $K$ denote the base field which we assume to be algebraically
closed of arbitrary characteristic. Given a distributive lattice
$\mathcal{L}$, let $X(\mathcal{L})$ denote the affine variety in
$\mathbb{A}^{\#\mathcal{L}}$ whose vanishing ideal is generated by
the binomials $X_\tau X_\varphi
-X_{\tau\vee\varphi}X_{\tau\w\varphi}$ in the polynomial algebra
$K[X_\alpha,\alpha\in\mathcal{L}]$ (here, ${\tau\vee\varphi}$
(resp. ${\tau\w\varphi}$) denotes the \emph{join} - the smallest
element of $\mathcal{L}$ greater than both $\tau,\varphi$ (resp.
the \emph{meet} - the largest element of $\mathcal{L}$ smaller
than both $\tau,\varphi$)). These varieties were extensively
studied by Hibi in  \cite{Hi} where Hibi proves that
$X(\mathcal{L})$ is a normal variety. On the other hand,
Eisenbud-Sturmfels show in \cite{ES} that a binomial prime ideal
is toric (here, ``toric ideal" is in the sense of \cite{St}). Thus
one obtains that $X(\mathcal{L})$ is a normal toric variety.  We
shall refer to such a $X\left(\mathcal{L}\right)$ as a \emph{Hibi
toric variety}.

For $\mathcal{L}$ being the Bruhat poset of Schubert varieties in
a minuscule $G/P$, it is shown in \cite{GLdef} that
$X(\mathcal{L})$ flatly deforms to ${\widehat{G/P}}$ (the cone
over $G/P$), i.e., there exists a flat family over
${\mathbb{A}}^1$ with ${\widehat{G/P}}$ as the generic fiber and
$X(\mathcal{L})$ as the special fiber. More generally, for a
Schubert variety $X(w)$ in a minuscule $G/P$, it is shown in
\cite{GLdef} that $X({\mathcal{L}}_w)$ flatly deforms to
${\widehat{X(w)}}$, the cone over $X(w)$ (here, ${\mathcal{L}}_w$
is the Bruhat poset of Schubert subvarieties of $X(w)$). In a
subsequent paper (cf. \cite{g-l}), the authors of loc.cit.,
studied the singularities of $X(\mathcal{L}), \mathcal{L}$ being
the Bruhat poset of Schubert varieties in the Grassmannian;
further, in loc.cit., the authors gave a conjecture (cf. \S 11 in
loc.cit; see also Remark \ref{conj11} of this paper for a
statement of this conjecture) giving a necessary and sufficient
condition for a point on $X(\mathcal{L})$ to be smooth, and proved
in loc.cit., the sufficiency part of the conjecture. Subsequently,
the necessary part of the conjecture was proved in \cite{Font} by
Batyrev et al. The toric varieties $X(\mathcal{L}), \mathcal{L}$
being the Bruhat poset of Schubert varieties in the Grassmannian
play an important role in the area of mirror symmetry; for more
details, see \cite{Font,Font2}. We refer to such an
$X(\mathcal{L})$ as a \emph{Grassmann-Hibi toric variety}
(abbreviated \emph{G-H toric variety}).

The proof (in \cite{g-l}) of the sufficiency part of the
conjecture of \cite{g-l} uses the Jacobian criterion for
smoothness, while the proof (in \cite{Font}) of the necessary part
of the conjecture of \cite{g-l} uses certain desingularization of
$X(\mathcal{L})$.

In this paper, our first main result is a description of the
singular locus of a G-H toric variety $X(\mathcal{L})$ in terms of
 the faces of the associated polyhedral cone; in
particular, we give a proof of the conjecture of \cite{g-l} using
just the combinatorics of the polyhedral cone associated to the
toric variety $X(\mathcal{L})$. Further, we prove (cf. Theorem
\ref{mainth}) that the singular locus of $X(\mathcal{L})$ is pure
of codimension 3 in $X(\mathcal{L})$, and that the generic
singularities are of cone type (more precisely, the singularity
type is the same as that at the vertex of the cone over the quadric
surface $x_1x_4-x_2x_3=0$ in $\mathbb{P}^3$); we also determine (in
loc.cit) the tangent cone to $\,X(\mathcal{L})$ at a maximal
singularity which turns out to be a toric variety (here, by a
``maximal singularity", we mean a point $P\in X(\mathcal{L})$ such
that the closure of the torus orbit through $P$ is an irreducible
component of Sing$\,X(\mathcal{L})$, the singular locus of
$X(\mathcal{L})$). We further obtain (cf. Corollary \ref{cor-8-3},
Theorem \ref{cat-mult}) an interpretation of the multiplicities at
some of the singularities as certain Catalan numbers in the case
of $\mathcal{L}$ being the Bruhat poset of Schubert varieties in
the Grassmannian of $2$-planes in $K^n$. We also present a product
formula (cf. Theorem \ref{cat-mult''}).

It turns out that the conjecture of \cite{g-l} does not extend to
a general $X(\mathcal{L})$ (see \S \ref{counter-example} for a
counter example). But, recently (cf.\cite{b-l}), we have proved
the conjecture of \cite{g-l} for other minuscule posets. It should
be remarked that this paper contains more results for the
Grassmann-Hibi toric varieties which can not be deduced from the
results of \cite{b-l}.

The singularities of the Hibi toric variety were also studied by
Wagner (cf. \cite{Wa}).

The sections are organized as follows: In \S \ref{first-sect}, we
recall generalities on toric varieties. In \S \ref{varlat}, we
recall generalities on distributive lattices. In \S \ref{s3}, we
introduce the Hibi toric variety $X({\L})$ and recall some results
from \cite{g-l,LM} on $X(\mathcal{L})$. In \S \ref{cone}, we
recall results from \cite{LM} on the polyhedral cone associated to
$X({\L})$. In \S \ref{sect-b-h}, we introduce the Grassmann-Hibi
toric variety. In \S \ref{singular-sect}, we prove our first main
result giving the description of Sing$\,X(\mathcal{L})$ in terms
of faces of the cone associated to $X(\mathcal{L})$;  and we
present our results on the tangent cones and deduce the
multiplicities at the associated points. In \S \ref{d=2-sect}, we
present our results for $X(\mathcal{L})$, $\mathcal{L}$ being the
Bruhat poset of Schubert varieties in the Grassmannian of
$2$-planes in $K^n$. In \S \ref{last-sect}, we present a formula
for the multiplicity at the unique $T$-fixed point of
$X(\mathcal{L})$, $\mathcal{L}$ being the Bruhat poset of Schubert
varieties in the Grassmannian of $d$-planes in $K^n$. In \S
\ref{conjectures}, we present a counter example to show that the
conjecture of \cite{g-l} does not extend to a general
$X(\mathcal{L})$; in this section, we also present two
conjectures.

\section{Generalities on toric varieties}\label{first-sect}

Since our main object of study is a certain affine toric variety,
we recall in this section some basic definitions on affine toric
varieties. Let $T=(K^*)^m$ be an $m$-dimensional torus.
\begin{defn}\label{equi} (cf. \cite{F}, \cite{KK})
An {\em equivariant affine embedding}  of a torus $T$ is an affine
variety $X\subseteq{\mathbb{A}}^l$ containing $T$ as an open
subset and equipped with a $T$-action $T\times X\to X$ extending
the action $T\times T\to T$ given by multiplication. If in
addition $X$ is normal, then $X$ is called an {\em affine toric
variety}.
\end{defn}

\subsection{The Cone Associated to a Toric Variety}\label{comb}
Let $M$ be the character group of $T$, and $N$ the $\Z$-dual of
$M$.  Recall (cf. \cite{F}, \cite{KK}) that there exists a
strongly convex rational polyhedral cone $\sigma \subset N_{\R}
(=N\otimes _{\Z}\R )$ such that
\[ K[X]=K[S_{\sigma}], \]
where $S_{\sigma}$ is the subsemigroup $\sigma^{\v}\cap M$,
$\sigma^{\v}$ being the cone in $M_{\R}$ dual to $\sigma$.  Note
that $S_{\sigma}$ is a finitely generated subsemigroup in $M$.

\subsection{Orbit Decomposition in Affine Toric Varieties}\label{orbit}
We shall denote $X$ also by $X_\sigma$.  We may suppose, without
loss of generality, that $\sigma$ spans $N_{\R}$ so that the
dimension of $\sigma$ equals $\dim N_{\R}=\dim T$. (Here, by
dimension of $\sigma$, one means the vector space dimension of the
span of $\sigma$.)

Let us first recall the definition of faces of a convex polyhedral
cone:
\begin{defn}
A \emph{face} $\tau$ of $\sigma$ is a convex polyhedral sub cone
of $\sigma$ of the form $\tau=\sigma\cap u^{\perp}$ for some $u\in
\sigma^{\v}$, and is denoted $\tau<\sigma$.  Note that $\sigma $
itself is considered a face.
\end{defn}

We have that $X_\tau$ is a principal open subset of $X_\sigma$,
namely, $$X_\tau=(X_\sigma)_u$$

\ni  Each face $\tau$ determines a (closed) point $P_\tau$ in
$X_\sigma$, namely, it is the point corresponding to the maximal
ideal in $K[X](=K[S_\sigma])$ given by the kernel of
$e_\tau:K[S_\sigma]\rightarrow K $, where for $u\in S_\sigma$, we
have
$$e_\tau(u)=\begin{cases}1,&{\mathrm{\ if\ }}u\in \tau^{\perp}\\
0,&{\mathrm{otherwise}} \end{cases}$$

\begin{rem}\label{coord} As a point in ${\mathbb{A}}^l$, $P_\tau$ may be identified with the
$l$-tuple with $1$ at the $i$-th place if $\chi_i$ is in
$\tau^{\perp}$, and $0$ otherwise (here, $\chi_i$ denotes the
weight of the $T$-weight vector $y_i $ - the class of $x_i$ in
$K[X_\sigma]$).
\end{rem}

\subsection{Orbit Decomposition}\label{orbit-decomp} Let $O_\tau$ denote
the $T$-orbit in $X_\sigma$ through $P_\tau$. We have the
following orbit decomposition in $X_\sigma$:
$$\begin{gathered}
X_\sigma=\bigcup_{\theta\le\sigma } O_\theta\\
{\overline{O_\tau}}=\bigcup_{\theta\ge\tau } O_\theta\\
\dim\,\tau + \dim\,O_\tau=\dim\,X_\sigma
\end{gathered}$$
See \cite{F}, \cite{KK} for details.

Thus $\tau\mapsto {\overline{O_\tau}}$ defines an order reversing
bijection between \{faces of $\sigma$\} and \{$T$-orbit closures
in $X_\sigma$\}. In particular, we have the following two extreme
cases:

\ni 1. If $\tau$ is the $0$-face, then $P_\tau=(1,\cdots,1)$ (as a
point in ${\mathbb{A}}^l$), and $O_\tau=T$, and is contained in
$X_\theta,$ for every $ \theta<\sigma$. It is a dense open orbit.

\ni 2. If $\tau=\sigma$, then $P_\tau =\left(0,\ldots ,0\right)$
(as a point in ${\mathbb{A}}^l$), and $O_\tau=\{P_\tau\}$, the
unique closed orbit.  (Note that since $\sigma$ spans $N_{\R}$,
$P_{\sigma}$ is a $T$-fixed point and is in fact the unique
$T$-fixed point in $X_\sigma$.)

For a face $\tau$, let us denote by $N_\tau$ the sublattice of $N$
generated by the lattice points of $\tau$. Let $N(\tau)=N/N_\tau$,
and $M(\tau)$, the ${{\mathbb{Z}}}$-dual of $N(\tau)$. For a face
$\theta$ of $\sigma$ such that $\theta$ contains $\tau$ as a face,
set
$$\theta_\tau:=(\theta+(N_\tau)_{{\mathbb{R}}})/(N_\tau)_{{\mathbb{R}}}.$$
Then the collection $\{\theta_\tau,\sigma\geq\theta\geq\tau\}$
gives the set of faces of the cone $\sigma_\tau (\subset
N(\tau)_{{\mathbb{R}}})$.

\begin{lem}[cf. \cite{F}]\label{closure} For a face $\tau<\sigma$, ${\overline{O_\tau}}$ gets
identified with the toric variety $Spec\,K[S_{\sigma_{\tau}}]$.
Further, $K[{\overline{O_\tau}}]=K[S_\sigma\cap\tau^{\perp} ]$.
\end{lem}

\section{Toric varieties associated to  finite distributive
lattices}\label{varlat} We shall now study a special class of
toric varieties, namely, the toric varieties associated to
distributive lattices. We shall first collect some definitions on
finite partially ordered sets. A partially ordered set is also
called a poset.
\begin{defn}
A finite poset $P$ is called\/ {\em bounded} if it has a unique
maximal, and a unique minimal element, denoted $\widehat{1}$ and
$\widehat{0}$ respectively.
\end{defn}
\begin{defn}
A totally ordered subset $C$ of a finite poset $P$ is called a\/
{\em chain}, and the number $\# C-1$ is called the\/ {\em length}
of the chain.
\end{defn}
\begin{defn}
A bounded poset $P$ is said to be\/ {\em graded} (or also {\em
ranked}) if all maximal chains have the same length (note that
$\widehat{1}$ and $\widehat{0}$ belong to any maximal chain).
\end{defn}
\begin{defn}
Let $P$ be a graded poset. The length of a maximal chain in $P$ is
called the\/ {\em rank} of $P$.
\end{defn}
\begin{defn}\label{5.5}
Let $P$ be a graded poset. For $\l,\m\in P$ with $\l\ge\m$, the
graded poset $\{\t\in P\mid\m\le\t\le\l\}$ is called the\/ {\em
interval from $\m$ to $\l$}, and denoted by $[\m,\l]$. The rank of
$[\m,\l]$ is denoted by $l_{\m}(\l)$; if $\m=\widehat{0}$, then we
denote $l_{\m}(\l)$ by just $l(\l)$.
\end{defn}

\begin{defn}
Let $P$ be a graded poset, and $\l,\m\in P$, with $\l\ge\m$. The
ordered  pair $(\l,\m)$ is called a\/ {\em cover} (and we also say
that $\l$\/ {\em  covers} $\m$) if $l_{\m}(\l)=1$.
\end{defn}

\subsection{Generalities on distributive lattices}\label{s2}
\begin{defn} A {\em lattice\/} is a partially ordered set $(\mathcal{L},\le)$ such that,
for every pair of elements $x,y\in \mathcal{L}$, there exist
elements $x\v y$ and $x\w y$, called the\/ {\em join},
respectively the\/ {\em meet} of $x$ and $y$, defined by:
\begin{gather}
x\v y\ge x,\ x\v y\ge y,\text{ and if  }z\ge x \text{ and }
z\ge y,\text{ then } z\ge x\v y,\notag\\
x\w y\le x,\ x\w y\le y,\text{ and if  }z\le x \text{ and } z\le
y,\text{ then } z\le x\w y.\notag
\end{gather}
\end{defn}
It is easy to check that the operations $\v$ and $\w$ are
commutative and associative.

\begin{defn}
Given a lattice ${\mathcal{L}}$, a subset ${\mathcal{L}}'\subset
{\mathcal{L}}$ is called a\/ {\em sublattice} of $\mathcal{L}$ if
$x,y\in{\mathcal{L}}'$ implies $x\w y\in{\mathcal{L}}'$, $x\v
y\in{\mathcal{L}}'$.
\end{defn}
\begin{defn}
Two lattices ${\mathcal{L}}_1$ and ${\mathcal{L}}_2$ are\/ {\em
isomorphic} if there exists a bijection
$\varphi:{\mathcal{L}}_1\to {\mathcal{L}}_2$ such that, for all
$x,y\in{\mathcal{L}}_1$,
\begin{equation*}
\varphi(x\v y)=\varphi(x)\v \varphi(y) \text{   and   }
\varphi(x\w y)=\varphi(x)\w \varphi(y).
\end{equation*}
\end{defn}
\begin{defn} A lattice is called\/ {\em distributive} if the following
identities hold:
\begin{align}
x\w (y\v z)&=(x\w y)\v (x\w z)\\
x\v (y\w z)&=(x\v y)\w (x\v z).
\end{align}
\end{defn}

\begin{defn}
An element $z$ of a lattice $\L$ is called\/ {\em
join-irreducible} (respectively\/ {\em meet-irreducible}) if
$z=x\v y$ (respectively $z=x\w y$) implies $z=x$ or $z=y$. The set
of join-irreducible (respectively meet-irreducible) elements of
$\L$ is denoted by $J\left(\L\right)$ (respectively
$M\left(\L\right)$), or just by $J$ (respectively $M$) if
 no confusion is possible.
\end{defn}

\begin{defn}
The set $J\left(\L\right)\cap M\left(\L\right)$ of join and
meet-irreducible elements is denoted by $JM\left(\L\right)$, or
just $JM$ if no confusion is possible.
\end{defn}
\begin{defn}
A subset $I$ of a poset $P$ is called an {\em ideal} of $P$ if for
all $x,\,y\in P$,
$$x\in I\text{ and }y\le x\text{ imply }y\in I.$$
\end{defn}
\begin{thm}[Birkhoff]\label{5.10}
Let $\L$ be a distributive lattice with $\hat{0}$, and $P$ the
poset of its nonzero join-irreducible elements. Then $\L$ is
isomorphic to the lattice of ideals of $P$, by means of the
lattice isomorphism
$$\a\mapsto I_{\a}=\{\t\in P\mid \t\le\a\},\qquad \a\in\L.$$
\end{thm}

\begin{defn}\label{diamond}
A quadruple of the form $(\t,\f,\t\v\f,\t\w\f)$, with $\t,\f\in\L$
non-comparable is called a\/ {\em diamond}, and is denoted by
$D(\t,\f,\t\v\f,\t\w\f)$ or also just $D(\t,\f)$. The unordered
pair $(\t,\f)$ (respectively $(\t\v\f,\t\w\f)$)) is called the
\emph{skew} (respectively \emph{main}) diagonal of the diamond
$D(\t,\f)$.
\end{defn}


The following Lemma is easily checked.
\begin{lem}
With the notations as above, we have

$(a)$ $J=\{\t\in\L\mid\text{there exists at most one cover of the
form } (\t,\l)\}$.

$(b)$ $M=\{\t\in\L\mid\text{there exists at most one cover of the
form } (\l,\t)\}$.
\end{lem}

For $\a\in\L$, let $I_\a$ be the ideal corresponding to $\a$ under
the isomorphism in Theorem \ref{5.10}; (we will also use the
notation $I\left(\alpha\right)$).
\begin{lem}[cf. \cite{LM}]\label{cover}
Let $(\t,\l)$ be a cover in ${\mathcal{L}}$. Then $I_{\t}$ equals
$I_{\l}\dot\cup\{\beta\}$ for some $\beta\in
J\left({\mathcal{L}}\right)$.
\end{lem}

\section{The variety $X(\L)$}\label{s3}

Consider the polynomial algebra $K[X_\alpha,\alpha\in\L]$; let
$I(\L)$ be the ideal generated by $\{X_\alpha X_\beta-X_{\alpha\v
\beta}X_{\alpha\w \beta}, \alpha,\beta\in\L\}$. Then one knows
(cf.\cite{Hi}) that $K[X_\alpha,\alpha\in\L]\,/I(\L)$ is a normal
domain; in particular, we have that $I(\L)$ is a  prime ideal. Let
$X(\L)$ be the affine variety of the zeroes in $K^l$ of $I(\L)$
(here, $l=\# \L$).
 Then $X(\L)$ is an affine normal variety defined by binomials; on
 the other hand, by \cite{ES}, a binomial prime ideal is a toric
 ideal (here, ``toric ideal" is in the sense of \cite{St}). Hence
  $X(\L)$ is a toric variety for the action by a suitable torus $T$.

In the sequel, we shall denote
$R\left(\mathcal{L}\right)=K[X_\alpha,\alpha\in\L]\,/I(\L)$.
Further, for $\alpha\in \L$, we shall denote the image of
$X_\alpha$ in $R\left(\mathcal{L}\right)$ by $x_{\alpha}$.

\begin{defn} The variety $X\left(\mathcal{L}\right)$ will be called
a \em{Hibi toric variety.}
\end{defn}

\begin{remark} An extensive study of $X\left(\mathcal{L}\right)$ appears first in \cite{Hi}.
\end{remark}

We have that $\dim X(\L)=\dim T$.

\begin{thm}[cf. \cite{LM}]\label{dim}
The dimension of $X(\L)$ is equal to  $\#J\left(\L\right)$.
\end{thm}

\begin{defn}
 For a finite distributive lattice $\L$, we call the
cardinality of $J\left(\L\right)$  the {\em dimension} of $\L$,
and we denote it by $\dim\L$. If $\L'$ is a sublattice of $\L$,
then the {\em codimension} of $\L'$ in $\L$ is defined as
$\dim\L-\dim\L'$.
\end{defn}

\begin{defn}\label{5.21}(cf. \cite{Wa})
A sublattice  $\L '$ of $\L$ is called an {\em embedded sublattice
of} $\L$ if
$$\t,\,\f\in\L,\quad\t\v\f,\,\t\w\f\in\L '\quad\Rightarrow\quad\t,\,\f\in\L '.$$
\end{defn}
Given a sublattice $\L '$ of $\L$, let us consider the variety
$X(\L ')$, and consider the canonical embedding $X(\L
')\hookrightarrow{\mathbb{A}}(\L
')\hookrightarrow{\mathbb{A}}(\L)$ (here ${\mathbb{A}}(\L
')={\mathbb{A}}^{\#\L '}$,
${\mathbb{A}}(\L)={\mathbb{A}}^{\#\L}$).

\begin{prop}[cf. \cite{g-l}]\label{5.22}
$X(\L ')$ is a subvariety of $X(\L)$ if and only if $\L '$ is an
embedded sublattice of $\L$.
\end{prop}

\subsection{Degeneration of $X(\L)$ to a monomial scheme}

\begin{defn}
A monomial $x_{\alpha_{1}}\cdots x_{\alpha_{r}}$ in $R({\L})$ of
degree $r$ is said to be \emph{standard} if
$\alpha_1\ge\cdots\ge\alpha_r$.
\end{defn}

\begin{thm}[cf.~\cite{g-l}]
Monomials standard on $X\left(\L\right)$ form a basis for
$R\left(\L\right)$.
\end{thm}

\begin{definition}\label{stanley-reisner} Let $S({\L})$ denote the \emph{Stanley-Reisner algebra} of ${\L}$,
namely,

\ni $K[X_\alpha,\alpha_i\in {\L}]/{\mathfrak{b}}\left({\L}\right)$
where ${\mathfrak{b}}\left({\L}\right)$ is the (monomial ideal)
ideal in $K[X_\alpha,\alpha_i\in {\L}]$ generated by
 $\{X_\alpha X_\beta,\alpha,\beta\mathrm{\ incomparable}\}$.
\end{definition}

Note that Spec$\,S({\L})$ is a monomial scheme. Using the standard
monomial basis, a flat deformation of $\Spec\,S\left({\L}\right)$
is constructed in \cite{g-l}:
\begin{thm}\label{flat}
There exists a flat family over ${\mathbb{A}}^{1}$ whose special
fiber ($t=0$) is  Spec$\,S({\L})$, and generic fiber ($t$
invertible) is $X\left({\L}\right)$.
\end{thm}

\subsection{Hilbert polynomial and Hilbert series:}\label{hilbert}
Recall the Hilbert function of a graded, finitely generated
$K$-algebra $R=\oplus_m\,R_m$;
$\phi_m(R):=\dim\,R_m$ is the {\emph{Hilbert function}} of $R$.
Recall that for $m\gg 0,\phi_m(R)$ is given by $P(m)$, where
$P(x)$ is a polynomial of degree equal to
$r:=\dim\,R-1(=\dim\,$Proj$\,R)$ and coefficients are in
${\mathbb{Q}}$; further, the leading coefficient of $P(x)$ is of
the form ${\frac{e_R}{r!}}$. The polynomial $P(x)$ is called
{\emph{the Hilbert polynomial}} of $R$, and $e_R$ is called the
{\emph{degree of $R$}} or also the {\emph{degree of $Proj\,R$}}.
The series $H_R(t):=\sum\,\phi_mt^m$ is called {\emph{the Hilbert
series of $R$}}.

\vs.1cm\ni{\textbf{Tangent cone and multiplicity:}} Let
$X=Spec\,B\hookrightarrow {\mathbb{A}}^l$ be an affine variety;
let $S$ be the polynomial algebra $K[X_1,\cdots,X_l]$. Let $P\in
X$, and let $M_P$ be the maximal ideal in $K[X]$ corresponding to
$P$ (we are concerned only with closed points of $X$). Let
$A={\mathcal{O}}_{X,P}$, the stalk at $P$; denote the unique
maximal ideal in $A$ by $M (=M_PB_{M_P})$. Then $Spec\,gr(A,M)$,
where
$gr(A,M)={\underset{j\in{\mathbb{Z}}_+}{\oplus}}\,M^j/M^{j+1}$ is
the {\emph{tangent cone}} to $X$ at $P$, and is denoted $TC_PX$.
Note that $$gr(B,M_P)=gr(A,M)$$

\ni The degree of the graded ring $gr(A,M)$ is defined as the
{\emph{multiplicity of }}$X$ at $P$, and is denoted mult\,$_PX$

\vs.1cm\ni\textbf{Fact:} With notation as above, if $B$ is graded,
then taking $P$ to be \textbf{0} (the origin), we have a natural
identification
\[{TC_{\mathbf{0}}X}\cong X\] Hence in this case, we have
\[\mathrm{mult}_{\mathbf{0}}X=e_B\]  See
\cite{Eisenbud} for details.

 \vs.2cm\ni\textbf{{Square-free monomial
schemes:}} Let $A=K[X_1,\cdots,X_n]$. Let $I$ be a square-free
monomial ideal in $A$; i.e., generated by a bunch of square-free
monomials ${\textbf{m}}_i,i\in J$, for some finite indexing set
$J$; let $R=A/I$. Let
\[I_i=\{j\mid X_j\mbox{ occurs in }\textbf{m}_i\}.\]

\begin{thm}[cf. \cite{LM}]\label{Hilbert} $\,$
\begin{enumerate}
\item $\phi_m(R)=\sum\,{{m-1}\choose{\#S-1}}$, the sum running
over all $S$ such that $S\not\supseteq I_i,\forall i$.  In
particular, the Hilbert polynomial of $R$ has degree $d-1$, where
$d=max\,\{\#S,S\not\supseteq I_i,\forall i\}$. \item
$H_R(t)=\sum\,{\frac{t^{\#S}}{(1-t)^{\#S}}}$, the sum running over
all $S$ such that $S\not\supseteq I_i,\forall i$. \item
$dim\,Spec\,R=d$. \item $deg\,R=\#\{S,S\not\supseteq I_i,\forall
i,\#S=d\}$
\end{enumerate}
\end{thm}

\begin{cor}
Let $P$ be a ranked poset, and $R(P)$ be the Stanley-Reisner
algebra of ${P}$ (cf. Definition~\ref{stanley-reisner}).  Then
$dim\,\Spec\,R(P)$ equals the cardinality of the set of elements
in a maximal chain in ${P}$, and $deg\,R\left(P\right)$ is equal
to the number of maximal chains in $P$.
\end{cor}

Now returning to our situation of $X\left({\L}\right)$, we get, in
view of Theorem \ref{flat}, and the above Corollary
\begin{thm}\label{DIM} $ $
\begin{enumerate}\item The dimension of $X\left({\L}\right)$ equals the cardinality of
the set of elements in a maximal chain in ${\L}$ \item The degree
of $K[X\left({\mathcal{L}}\right)]$ is equal to the number of
maximal chains in $\mathcal{L}$.
\end{enumerate}
\end{thm}
Combining the above Theorem with Theorem \ref{dim}, and Lemma
\ref{cover}, we obtain
\begin{thm}\label{card}$\,$
\begin{enumerate}\item The cardinality of $\,J\left(\L\right)$ equals the cardinality of the set of elements in a
maximal chain in ${\L}$. \item Fix any maximal chain
$\beta_1<\cdots<\beta_d$ in $X({\L})$, $d$ being $\#
J\left(\L\right)$. Let $\gamma_{i+1}$ be the element of
$J\left(\L\right)$ corresponding to the cover
$(\beta_{i+1},\beta_i)$, $i\ge 1$ (cf. Lemma~\ref{cover}); set
$\gamma_1=\beta_1$. Then
$J\left(\L\right)=\{\gamma_1,\cdots,\gamma_d\}$
\end{enumerate}
\end{thm}
Now $X({\L})$ being of cone type, i.e, the vanishing ideal of
$X({\L})$ (as a subvariety of ${\mathbb{A}}^{\#{\L}}$) is
homogeneous, we have (in view of \S \ref{hilbert}, Fact) the
following
\begin{thm}\label{DIM'} $ $
\begin{enumerate}\item  We have an identification of $X({\L})$ with
$TC_{P_{\sigma}}X({\L})$ (the tangent cone to $X({\L})$ at
$P_\sigma$, origin being identified with $P_\sigma$)\item
$mult_{P_{\sigma}}X({\L})$ equals the degree of
$K[X\left({\mathcal{L}}\right)]$ ( $=\# \{\mathrm{maximal\ chains\
in}\ {\mathcal{L}}\}$).
\end{enumerate}
\end{thm}
\section{Cone and dual cone of $X\left({\mathcal{L}}\right)$:}\label{cone}
Denote the poset of join-irreducibles by $J\left({\L}\right)$.
Denote by ${\mathcal{I}}$ the poset of ideals of
$J\left({\L}\right)$. For  $A\in {\mathcal{I}}$, denote by
${\mathbf{m}}_A$ the monomial:
$${\mathbf{m}}_A:={\underset{\tau\in A}{\prod}} \,y_\tau$$ in the
polynomial algebra $K[y_\tau,\tau\in J\left({\L}\right)]$. If
$\alpha$ is the element of $\L$ such that $I_\alpha=A$ (cf.
Theorem \ref{5.10}), then we shall denote ${\mathbf{m}}_A$ also by
${\mathbf{m}}_\alpha$. Consider the surjective algebra map
$$F:K[X_\alpha,\alpha\in{\mathcal{L}}]\rightarrow
K[{\mathbf{m}}_A,A\in{\mathcal{I}}],\,X_\alpha\mapsto
{\mathbf{m}}_A,\ A=I_\alpha
$$

\begin{thm}[cf. \cite{Hi}, \cite{LM}]\label{main1} We have an isomorphism
$K[X\left({\L}\right)]\cong K[{\mathbf{m}}_A,A\in{\mathcal{I}}]$.
\end{thm}

Let us denote the torus acting on the toric variety
$X\left({\L}\right)$ by $T$; note that $\dim\,T=\#
J\left({\L}\right)=d$, say. Identifying $T$ with $(K^*)^d$, let
$\{f_z,z\in J\left({\L}\right)\}$ denote the standard
${\mathbb{Z}}$-basis  for $X(T)$, namely, for $t=(t_z, z \in
J\left({\L}\right))$, $f_z(t)=t_z$. Denote $M:=X(T)$; let $N$ be
the ${\mathbb{Z}}$-dual of $M$, and $\{e_y,y\in
J\left({\L}\right)\}$ be the basis of $N$ dual to $\{f_z,z\in
J\left({\L}\right)\}$. For $A\in {\mathcal{I}}$, set
$$f_A:={\underset{z\in A}{\sum}} \,f_z$$ Let
$V=N_{{\mathbb{R}}}(=N\otimes_{{\mathbb{Z}}} {\mathbb{R}})$. Let
$\sigma \subset V$ be the cone such that
$X\left({\mathcal{L}}\right)=X_\sigma$. Let $\sigma^{\v}\subset
V^*$ be the cone dual to $\sigma$. Let $S_\sigma=\sigma^{\v}\cap
M$, so that $K[X\left({\mathcal{L}}\right)]$ equals the semi group
algebra $K[S_\sigma]$.

As an immediate consequence of Theorem \ref{main1}, we have
\begin{prop}\label{semi}
The semigroup $S_\sigma$ is generated by $f_A, A\in
{\mathcal{I}}$.
\end{prop}

Let $M(J\left({\L}\right))$ be the set of maximal elements in the
poset $J\left({\L}\right)$. Let $Z(J\left({\L}\right))$ denote the
set of all covers in the poset $J\left({\L}\right)$ (i.e.,
$(z,z'), z>z'$ in the poset $J\left({\L}\right)$, and there is no
other element $y\in J\left({\L}\right)$ such that $z>y>z'$). For a
cover $(y,y')\in Z(J\left({\L}\right))$, denote
$$v_{y,y'}:=e_{y'}-e_{y}$$

\begin{prop}[cf. \cite{LM}]\label{prop-cpc} The cone $\sigma$ is
generated by $\{e_{z},z\in M(J\left({\L}\right)),\,v_{y,y'},
(y,y')\in Z(J\left({\L}\right))\}$.
\end{prop}

\subsection{Analysis of faces of $\sigma$}\label{anal}
We shall concern ourselves
 just with the closed points in $X\left( {\L}\right)$.
So in the sequel, by a point in $X\left( {\L}\right)$, we shall
mean a closed point. Let $\tau$ be a face of $\sigma$. Let
$P_{\tau}$ be the distinguished point of $O_{\tau}$ with the
associated maximal ideal being the kernel of the map
$$\begin{gathered}K[S_\sigma]\rightarrow K,\\
u\in S_\sigma,u\mapsto\begin{cases} 1,&\mathrm{\ if\ }
u\in\tau^{\perp}\\
0,&\mathrm{\ otherwise}
\end{cases}
\end{gathered}$$

Then for a point $P\in X\left( {\L}\right)$ (identified with a
point in ${\mathbb{A}}^l$, $l=\#\L$), denoting by $P(\alpha)$, the
$\alpha$-th co-ordinate of $P$, we have,
$$P_{\tau}(\alpha)=\begin{cases} 1,&\mathrm{\ if\ }
f_{I_{\alpha}}\in\tau^{\perp}\\
0,&\mathrm{\ otherwise}
\end{cases}$$
With notation as above, let
$$D_{\tau}=\{\alpha\in{\L}\,|\,P_{\tau}(\alpha)\neq 0\}$$ We have,

\subsection{The bijection $\mathcal{D}$}\label{bij} (cf. \cite{LM})
We have a bijection
$$\mathcal{D}:\{\mathrm{\ faces\ of\
}\sigma\}\,{\buildrel{bij}\over{\leftrightarrow}}\, \{\mathrm{\
embedded\ sublattices \ of\ }{\L}\},\ \mathcal{D}(\tau)=D_\tau$$
\begin{prop}\label{orbitdim}[cf. \cite{LM}]
Let $\tau$ be a face of $\sigma$. Then we have
${\overline{O_{\tau}}}=X\left({D_{\tau}}\right)$.
\end{prop}


\section{The Distributive Lattice $I_{d,n}$ and the
Grassmann-Hibi \\ Toric Variety}\label{sect-b-h}

We now turn our focus to
\[{\mathcal{L}}=I_{d,n}=\{x=(i_1,\ldots i_d)\mid 1\leq i_1<\ldots <i_d\leq n\}.\]
The partial order $\ge$ on $I_{d,n}$ is given by
$$(i_1,\dots,i_d)\ge (j_1,\dots,j_d)\iff i_1\ge j_1,\dots,i_d\ge j_d.$$
For $x\in I_{d,n}$, we denote the $j$-th entry in $x$ by $x(j)$,
$1\le j\le d$.

\begin{lemma} $I_{d,n}$ is a distributive lattice.
\end{lemma}
\begin{proof} To see that $I_{d,n}$ is a lattice, we simply need to see that there is a well defined meet and join.
We see that
\[\left(i_1,\ldots i_d\right)\v\left(j_1,\ldots j_d\right)=\left(M_1,\ldots M_d\right),\]
where $M_k=$max$\{i_k,j_k\}$.  It's clear that this is an element
of the lattice.  Similarly,
\[\left(i_1,\ldots i_d\right)\w\left(j_1,\ldots j_d\right)=\left(m_1,\ldots m_d\right),\]
for $m_k=$min$\{i_k,j_k\}$.  For proof of distributivity, see
\cite{Hiller}.
\end{proof}

\begin{remark} It is a well known fact (see \cite{L-G} for example) that the partially ordered set $I_{d,n}$ is isomorphic to the poset determined by the set of Schubert varieties in $G_{d,n}$, the Grassmannian of $d$-dimensional subspaces in an $n$-dimensional space, the Schubert varieties being partially ordered by inclusion.  \end{remark}

In the lemma below, by a\/ {\em segment} we shall mean a set
consisting of consecutive integers.
\begin{lem}[cf. \cite{g-l}]\label{8.2}
We have

$(a)$ The element $\t=(i_1,\dots,i_d)$ is join-irreducible if and
only if either $\t$ is a segment, or $\t$ consists of two disjoint
segments $(\m,\n)$, with $\m$ starting with $1$.

$(b)$ The element $\t=(i_1,\dots,i_d)$ is meet-irreducible if and
only if either $\t$ is a segment, or $\t$ consists of two disjoint
segments $(\m,\n)$, with $\n$ ending with $n$.

$(c)$ The element $\t=(i_1,\dots,i_d)$ is join-irreducible and
meet-irreducible if and only if either $\t$ is a segment, or $\t$
consists of two disjoint segments $(\m,\n)$, with $\m$ starting
with $1$ and $\n$ ending with $n$.
\end{lem}

\begin{defn} We shall denote $X\left(I_{d,n}\right)$ by just
$X_{d,n}$ and will refer to it as a \emph{Grassmann-Hibi toric
variety}, or a G-H toric variety for short.
\end{defn}

\section{Singular Faces of the G-H toric variety $X_{d,n}$}\label{singular-sect}

Let ${\L}$ represent the distributive lattice $I_{d,n}$. Let
$J\left({\L}\right)$ be the join irreducible elements of ${\L}$.
  From Lemma~\ref{8.2}, we have that the elements of
   $J\left({\L}\right)$ are of two types, elements of Type I
   consist of segments of ${\L}$, and those of Type II consist
   of all elements $\tau\in {\L}$
   such that $\tau =\left(\mu ,\nu\right)$ where $\mu$ is a
   segment beginning at $1$ and $\nu$ is just a segment.

Since the generators of the cone $\sigma$ are determined by
$J\left({\L}\right)$ (cf. Proposition~\ref{prop-cpc}), we will
often consider $J\left({\L}\right)$ as a partially ordered set
with the partial order induced from $\L$. Notice that
$J\left({\L}\right)$ has one maximal element, which is also the
maximal element of ${\L}$: $\hat{1}=\left(n-d+1,\ldots n\right)$;
and $J\left(\L\right)$ has one minimal element, which is also the
minimal element of ${\L}$: $\hat{0}=(1,\ldots d)$. Viewing
$J\left({\L}\right)$ as a poset, for each element $x$ of
$J\left({\L}\right)$, there are at most two covers of the form
$\left(y,x\right)$.

For example, if $x=\left(1,\ldots k, l+1,\ldots l+d-k\right)\in
J\left({\L}\right)$, we have $y=(1,\ldots k,l+2,\ldots l+d-k+1)$
and $y'=(1,\ldots k-1,l,\ldots l+d-k)$ forming the two covers of
$x$ in $J\left({\L}\right)$ (if $k=1$, then $y'=(l,\ldots
l+d-1)$). If $l=n-d+k$, or if $x$ is of Type I, then $x$ has only
one cover.

\begin{lemma}[cf. \cite{Proctor}] The partially ordered set $J\left({\L}\right)$ is a distributive lattice.
\end{lemma}

\begin{remark}\label{diagram-remark} As a lattice, $J\left({\L}\right)$ looks like a tessellation of diamonds in the shape of a rectangle, with sides of length $d-1$ and $n-d-1$.  For example, let $d=3$ and $n=7$.  Then $J\left({\L}\right)$ is the lattice \small
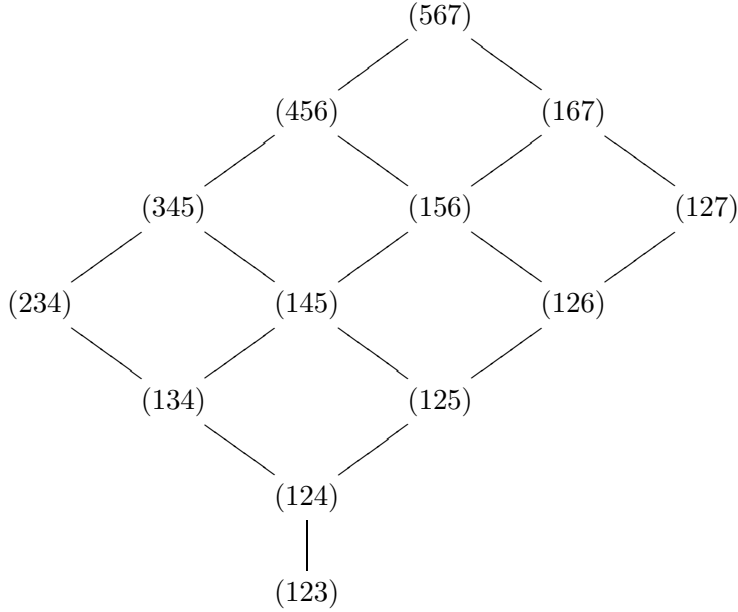
\[ \xymatrix@-5pt{
 & & & (567) & &\\
 & & (456)\ar@{-}[ur] & & (167)\ar@{-}[ul] & \\
 &(345)\ar@{-}[ur] & & (156)\ar@{-}[ul]\ar@{-}[ur] & & (127)\ar@{-}[ul]\\
 (234)\ar@{-}[ur]& & (145)\ar@{-}[ul]\ar@{-}[ur] & & (126)\ar@{-}[ul]\ar@{-}[ur] & \\
 & (134)\ar@{-}[ul]\ar@{-}[ur] & & (125)\ar@{-}[ul]\ar@{-}[ur]& & \\
 & & (124)\ar@{-}[ur]\ar@{-}[ul]& & & \\
 & & (123)\ar@{-}[u] & & &
 } \] \normalsize
\end{remark}

\noindent\textbf{The Face $\sigma_{ij}$:} As in \S \ref{cone}, let
$\sigma$ be the cone associated to $X\left(\L\right)$.

\begin{definition} For $1\leq i\leq n-d-1$, $1\leq j\leq d-1$, let \[\mu_{ij}=\left(1,\ldots j,\,i+j+1,\ldots i+d\right),\mbox{ and}\]
\[\lambda _{ij}=\left(i+1,\ldots i+j,n+1+j-d,\ldots n\right).\]
Define
\[\mathcal{L}_{ij}=\mathcal{L}\,\backslash\, [\mu_{ij},\lambda_{ij}].\]
\end{definition}

\begin{remark} (i) By \cite{g-l}, Lemma 11.5, we have that
${\mathcal{L}}_{ij}$ is an embedded sublattice.

\ni (ii) For $\alpha$, $\beta\in J\left({\L}\right)$,
$\alpha\wedge\beta=\mu_{ij}$ for some $1\leq i \leq n-d-1,\, 1\leq
j\leq d-1$; thus every diamond in $J\left({\L}\right)$ has a
$\mu_{ij}$ as its minimal element.
\end{remark}

\begin{definition} Let $\sigma_{ij}$ be the face of $\sigma$ for which $D_{\sigma_{ij}}=\Lij$.  \end{definition}

\begin{definition}\label{def-sing-face} A face $\tau$ of $\sigma$ is a \emph{singular} (resp. \emph{non-singular}) face if $P_\tau$ is a singular (resp. non-singular) point of $X_{\sigma}$.
\end{definition}
Our first result is that $\sigma_{ij}$ is a singular face.  To
prove this, we first determine a set of generators for
$\sigma_{ij}$.
\begin{definition}\label{gen} Let us denote by $W(\sigma)$ (or just $W$) the set of generators for $\sigma$ as
described in Proposition \ref{prop-cpc}.  For a face $\tau$ of
$\sigma$, define
\[W\left(\tau\right)=\{v\in W\mid f_{I\left(\alpha\right)}
\left(v\right)=0,\,\forall\,\alpha\in D_\tau\}.\] (Here, $D_\tau$
is as in \S \ref{bij} ). Then $W\left(\tau\right)$ gives a set of
generators for $\tau$.
\end{definition}
\subsection{Determination of $W\left(\sigma_{ij}\right)$}\label{Wij} It will aid our proof below to notice a few facts about the generators of $\sigma_{ij}$.  First of all, $e_{\hat{1}}$ is not a generator for any $\sigma_{ij}$; because $\hat{1}\in \Lij$ for all $1\leq i\leq n-d-1$, $1\leq j\leq d-1$, and $e_{\hat{1}}$ is non-zero on $f_{I_{\hat{1}}}$.  Similarly, for the cover $(y',\hat{0})$, where $y'=(1,\ldots d-1, d+1)$; $e_{\hat{o}}-e_{y'}$ is not a generator for any $\sigma_{ij}$.

Secondly, for any cover $(y',y)$ in $J\left({\L}\right)$,
$e_y-e_{y'}$ is not a generator of $\sigma_{ij}$ if
$y\in\mathcal{L}_{ij}$, because
$f_{I\left(y\right)}\left(e_y-e_{y'}\right)\neq 0$.  Thus, in
determining elements of $W\left(\sigma_{ij}\right)$, we need only
be concerned with elements $e_y-e_{y'}$ such that $y\in
J\left(\L\right)\cap [\mu_{ij},\lambda_{ij}]$.  The elements of
$J\left(\L\right)\cap [\mu_{ij},\lambda_{ij}]$ are
\begin{eqnarray*}
 & & y_t=(1,\ldots ,j,i+j+1+t,\ldots ,i+d+t)\mbox{ for } 0\leq t\leq n-d-i\\
 & & z_t=(1,\ldots ,j-t,i+j+1-t,\ldots ,i+d)\mbox{ for } 0\leq t\leq j
\end{eqnarray*}
Note that $y_0=z_0=\mu_{ij},z_j=(i+1,\ldots ,i+d)$. We shall now
prove in the Theorem below that $W\left(\sigma_{ij}\right)$
consists of precisely four elements, forming a diamond in the
distributive lattice $J(\mathcal{L})$ with $\mu_{ij}$ as the
smallest element.

\begin{theorem}\label{main} $W\left(\sigma_{ij}\right)=\{ e_{\mu_{ij}}-e_A,\,e_{\mu_{ij}}-e_B,\,e_A-e_C,\,e_B-e_C\}$, where $A$, $B$, and $C$ are defined in the proof.
\end{theorem}
\begin{proof} We  divide the proof into two cases:
$j=1$, and $j>1$.\\

\noindent \textbf{Case 1:} Let $j=1$ and $1\leq i\leq n-d-1$. Here
we have
\[\mu_{ij}=\left(1,i+2,\ldots , i+d\right)\mbox{ and }\lambda_{ij}=\left(i+1,n-d+2,\ldots , n\right).\]
As discussed previously, we find that $\mu_{ij}$ is covered in
$J\left({\L}\right)$ by $A=(i+1,\ldots , i+d)$ and $B=(1,i+3,\ldots ,
i+d+1)$.  We have that both $A$ and $B$ are in the interval
$[\mu_{ij},\lambda_{ij}]$. Let $C$ be the join of $A$ and $B$ in
the lattice $J\left({\L}\right)$, i.e. \[C=\left(i+2,\ldots ,
i+d+1\right).\] Note that $\left(C,A\right)$ and
$\left(C,B\right)$ are covers in
$J\left({\L}\right)$.    \\

We first observe the following:
$$\mathrm{If\ }x=(x_1,\cdots,x_d)\in\mathcal{L}_{ij}\mathrm{\ is\ }\ge\mu_{ij},\mathrm{\ then\ }
x\ge C \leqno{(*)}$$ (this follows since the facts that
$x\ge\mu_{ij}$ and $x\in\mathcal{L}_{ij}$ imply that
$x\not\le\lambda_{ij}$, and hence $x_1\ge i+2$).

 \noindent\textbf{Claim (i):}
$e_{\mu_{ij}}-e_A$ and $e_{\mu_{ij}}-e_B$ are both in
$W\left(\sigma _{ij}\right)$.

We shall prove the claim for $e_{\mu_{ij}}-e_A$ (the proof for
$e_{\mu_{ij}-e_B}$ is similar).  To prove that $e_{\mu_{ij}}-e_A$
is in $W\left(\sigma_{ij}\right)$, we need to show that there does
not exist $x=\left(x_1,\ldots , x_d\right)\in\Lij$ such that $x\geq
\mu_{ij}$ and $x\not\geq A$.  But this follows from $(\ast )$ (note that
$(\ast )$ implies that if
$x=(x_1,\cdots,x_d)\in\mathcal{L}_{ij}\mathrm{\ is\
}\ge\mu_{ij},\mathrm{\ then\ } x\ge A$).

\noindent\textbf{Claim (ii):} $e_A-e_C$ and $e_B-e_C$ are in
$W\left(\sigma_{ij}\right)$.

The proof is similar to that of Claim (i).  Again, we shall prove
the result for $e_A-e_C$, (the proof for $e_B-e_C$ is similar). We
need to show that there does not exist $x=(x_1,\ldots
,x_d)\in\Lij$ such that $x\geq A$, but $x\not\geq C$. Again this
follows from $(\ast )$ (note that $x\ge A$ implies in particular that
$x\ge\mu_{ij}$).

 \noindent\textbf{Claim (iii):}
$W\left(\sigma_{ij}\right)=\{
e_{\mu_{ij}}-e_A,\,e_{\mu_{ij}}-e_B,\,e_A-e_C,\,e_B-e_C\}$.

  In the case under consideration, $j$ being 1, the only elements of $J\left({\L}\right)\cap [\mu_{ij},\lambda_{ij}]$ are of the form
\[y_t=\left(1,i+t+2,\ldots , i+d+t\right)\mbox{ for }0\leq t\leq n-d-i;\]
\[\mbox{ and }z_1=(i+1,\ldots , i+d).\]

Let $y'_t=\left(i+t+1,\ldots , i+d+t\right)$ for $1\leq t\leq
n-d-i$; thus we have covers of the type $\left(y'_t,y_t\right)$
for $0\leq t\leq n-d-i$, and $\left(y_{t+1},y_t\right)$ for $0\leq
t\leq n-d-1-i$.  Note that $y_0=\mu_{ij}$, $y_1=B$, $z_1=y_0'=A$
and $y_1'=C$.  In Claims (i) and (ii), we have shown that the
covers $(y_1,y_0)$, $(y'_0,y_0)$, $(y'_1,y_1)$, and $(y_1',z_1)$
yield elements of $W\left(\sigma_{ij}\right)$.  Note also that $C$
is the only cover of $A$.  Hence, it only remains to show that
$e_{y_t}-e_{y_t'}\not\in W\left(\sigma_{ij}\right)$ for $2\leq
t\leq n-d-i$, and $e_{y_t}-e_{y_{t+1}}\not\in
W\left(\sigma_{ij}\right)$ for $1\leq t\leq n-d-1-i$. For each of
these covers, we shall exhibit a $x\in\mathcal{L}_{ij}$ such that
$f_{I(x)}$ is non-zero on the cover under consideration.

Define $x_t=\left(i+t,i+t+2,\ldots , i+d+t\right)$; we have that $x_t\in\Lij$ for $2\leq t\leq n-d-i$.  Further, $f_{I\left(x_t\right)}$ is non-zero on $e_{y_{t}}-e_{y_{t+1}}$ for $2\leq t\leq n-d-i-1$ and $e_{y_t}-e_{y'_t}$ for $2\leq t\leq n-d-i$.  For $(y_2,y_1)$, note that $C\in\Lij$, and $f_{I\left(C\right)}$ is non-zero on $e_{y_1}-e_{y_2}$.  This completes the proof of Case 1.\\

\noindent\textbf{Case 2:}  Now let $2\leq j\leq d-1$, $1\leq i\leq
n-d-1$.  We have
\[\mu_{ij}=\left(1,\ldots , j,i+j+1,\ldots , i+d\right),\mbox{ and}\]
\[\lambda _{ij}=\left(i+1,\ldots , i+j,n+1+j-d,\ldots , n\right).\]

As in Case 1, we will look for covers of $\mu_{ij}$ in
$J\left({\L}\right)$.  They are $A=(1,\ldots , j-1,i+j,\ldots , i+d)$,
and $B=(1,\ldots , j,i+j+2,\ldots , i+d+1)$.  Define $C$ to be the
join of $A$ and $B$ in the lattice $J\left({\L}\right)$, thus
\[C=\left(1,\ldots , j-1,i+j+1,\ldots , i+d+1\right).\]

\noindent\textbf{Claim (iv):}
$\{e_{\mu_{ij}}-e_A,\,e_{\mu_{ij}}-e_B,\,e_A-e_C,\,e_B-e_C\}$ are
in $W\left(\sigma_{ij}\right)$.

We first observe the following:
$$\mathrm{If\ }x=(x_1,\cdots,x_d)\in\mathcal{L}_{ij}\mathrm{\ is\ }
\ge\mu_{ij},\mathrm{\ then\ } x\ge C \leqno{(**)}$$ For, suppose
that $x\not\ge C$. Now the facts that $x\ge\mu_{ij}$ and
$x\in\mathcal{L}_{ij}$ imply that $x\not\le\lambda_{ij}$, and thus
$x_l>i+l$ for some $1\leq l\leq j$.  Also, $x\not\geq C$; hence
$x_k <i+k+1$ for some $j\leq k\leq d$. Therefore
\[x=(x_1,\ldots , x_{l-1},\, x_l>i+l,\, x_{l+1}>i+l+1,\ldots
,x_{k-1}>i+k-1,\]\[ i+k+1>x_k>i+k,\ldots )\]  Clearly, no such
$x_k$ exists, and thus $(**)$ follows.

By $(**)$, we have that if $x\in \Lij$ is such that $x\geq\mu_{ij}$, then $x\geq A,\,B,$ and $C$.  Hence Claim (iv) follows.\\

\noindent\textbf{Claim (v):}
$W\left(\sigma_{ij}\right)=\{e_{\mu_{ij}}-e_A,\,e_{\mu_{ij}}-e_B,\,e_A-e_C,\,e_B-e_C\}$.

As in Claim (iii), we will show that all other covers in
$J\left({\L}\right)$ of the form $\left(y',y\right)$ where $y\in
J\left({\L}\right)\cap [\mu_{ij},\lambda_{ij}]$ are not in
$W\left(\sigma_{ij}\right)$.  As in \S \ref{Wij}, all of the
elements of $J\left({\L}\right)\cap [\mu_{ij},\lambda_{ij}]$ are
\begin{eqnarray*}
 & & y_t=(1,\ldots ,j,i+j+1+t,\ldots ,i+d+t)\mbox{ for } 0\leq t\leq n-d-i\\
 & & z_t=(1,\ldots ,j-t,i+j+1-t,\ldots ,i+d)\mbox{ for } 0\leq t\leq j
\end{eqnarray*}
(note that $z_j=(i+1,\ldots ,i+d)$). We will examine covers of
these elements; note that $y_0=z_0=\mu_{ij}$, $z_1=A$, and
$y_1=B$.

Let $z_t'=(1,\ldots ,j-t,i+j+2-t,\ldots ,i+d+1)$ for $1\leq t\leq
n-d-i$, and $z_j'=(i+2,\ldots ,i+d+1)$.  First, we want to show
that the covers $(z_{t+1},z_t)_{1\leq t\leq j-1}$ and
$(z_t',z_t)_{2\leq t\leq j}$ do not yield elements in
$W\left(\sigma_{ij}\right)$.  Note that
$\left(z_1',z_1\right)=\left(C,A\right)$, and $e_A-e_C\in
W\left(\sigma_{ij}\right)$.  Also, $C\in\Lij$, and
$f_{I\left(C\right)}\left(e_{z_1}-e_{z_2}\right)$ is non-zero;
thus $e_{z_1}-e_{z_2}\not\in W\left(\sigma_{ij}\right)$, and we
may restrict our concern to $t\geq 2$.  Let
\[x_t=(1,\ldots , j-t,i+j+1-t,n-d+j-t+2,\ldots , n) \mbox{ for }2\leq t\leq j-1,\]
\[x_j=(i+1,n-d+2,\ldots , n).\]
Now, on the interval $2\leq t\leq j-1$, we have the following
facts:
\begin{enumerate}
\item $x_t\geq z_t$, \item $x_t\not\geq z_{t+1}$, \item
$x_t\not\geq z_t'$. \item $x_t\not\leq \lambda_{ij}$
\end{enumerate}
Facts (1), (3), and (4) above hold for the case $t=j$, it is just
a separate check.  Therefore, for $2\leq t\leq j$ (resp. $2\leq
t\leq j-1$), we have $x_t\in\Lij$ and $f_{I\left(x_t\right)}$ is
non-zero on $e_{z_t}-e_{z_t'}$ (resp. $e_{z_{t}}-e_{z_{t+1}}$).

Next, we must concern ourselves with covers involving $y_t$.
Define
\[y'_t=\left(1,\ldots , j-1,i+j+t,\ldots , i+d+t\right)\mbox{ for }1\leq t\leq n-d-i.\]
To complete Claim (v), we must show that the covers
\[\left(y_{t+1},y_t\right)_{1\leq t\leq n-d-i-1}\mbox{ and }
(y'_t,y_t)_{2\leq t\leq n-d-i}\]
do not yield elements of
$W\left(\sigma _{ij}\right)$. Note that
$\left(y'_1,y_1\right)=\left(C,B\right)$ and thus does yield an
element of $W\left(\sigma_{ij}\right)$. Also,
$f_{I\left(C\right)}\left(e_{y_1}-e_{y_2}\right)$ is non-zero,
 therefore we can restrict our attention to $t\geq 2$.
 Let $x'_t=(1,\ldots , j-1,i+j+1,i+j+t+1,\ldots , i+d+t)$.
 On the interval $2\leq t\leq n-d-i$, we have the following facts:
\begin{enumerate}
\item $x'_t\geq y_t$, \item $x'_t\not\geq y'_t$, \item
$x'_t\not\geq y_{t+1}$, (for $t\leq n-d-i-1$), \item $x'_t\not\leq
\lambda_{ij}$.
\end{enumerate}
Therefore, on the interval $2\leq t\leq n-d-i$ (resp. $2\leq t\leq
n-d-i-1$), we have $x'_t\in\Lij$ and $f_{I\left(x'_t\right)}$ is
non-zero on $e_{y_t}-e_{y'_t}$ (resp. $e_{y_{t}}-e_{y_{t+1}}$).

This completes Claim (v), Case 2, and the proof of the theorem.
\end{proof}

\begin{remark}\label{sigma-ij-diamond} The face $\sigma_{ij}$ corresponds to the diamond in $J\left({\L}\right)$:
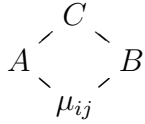
\[\xymatrix@=2.5pt{
 & C & \\
A\ar@{-}[ur] &  & B\ar@{-}[ul]\\
 & \mu_{ij}\ar@{-}[ul]\ar@{-}[ur] &
 } \]
which is a poset of rank 2.
\end{remark}

\begin{lemma}\label{dim3} $\dim\left(\sigma_{ij}\right)=3$.\end{lemma}
\begin{proof}  We have a set of generators for $\sigma_{ij}$, namely $\{e_{\mu_{ij}}-e_A,\,e_{\mu_{ij}}-e_B,\,e_A-e_C,\,e_B-e_C\}$.  We can see that a subset of three of these generators is linearly independent.  Thus if the fourth generator can be put in terms of the first three, the result follows.  Notice that
\[\left(e_{\mu_{ij}}-e_A\right)-\left(e_{\mu_{ij}}-e_B\right)+\left(e_A-e_C\right)=e_B-e_C .\]
\end{proof}

As an immediate consequence of Theorem \ref{main} and Lemma
\ref{dim3}, we have the following
\begin{theorem}\label{quadric} We have an identification of the (open)
affine piece in $X(\L)$ corresponding to the face $\sigma_{i,j}$
with the product $Z\times (K^*)^{\# J({\L})-3}$, where $Z$ is the
cone over the quadric surface $x_1x_4-x_2x_3=0$ in $\mathbb{P}^3$.
\end{theorem}

We now prove two lemmas which hold for a general toric variety.

\begin{lemma}\label{l-1-tau} Let $X_{\tau}$ be an affine toric variety with $\tau$ as the associated cone.  Then $X_{\tau}$ is a non-singular variety if and only if it is non-singular at the distinguished point $P_{\tau}$.
\end{lemma}
\begin{proof} Only the implication $\Leftarrow$ requires a proof.
Let then $P_\tau$ be a smooth point.  Let us assume (if possible)
that Sing$\,X_\tau\neq\varnothing$.  We have the following facts:

$\bullet$ Sing$\,X_\tau$ is a closed $T$-stable subset of
$X_\tau$.

$\bullet$ $P_\tau\in \overline{O_\theta}$, for every face $\theta$
of $\tau$ (see \S \ref{orbit-decomp}); in particular,
$P_\tau\in\overline{O_\theta}$, for some face $\theta$ such that
$P_\theta$ is a singular point, (such a $\theta$ exists, since by
our assumption Sing $X_\tau$ is non-empty). Therefore we obtain
that $P_\tau\in $ Sing $X_\tau$, a contradiction. Hence our
assumption is wrong and the result follows.
\end{proof}

\begin{lemma}\label{basis-sing} Let $\tau$ be a face of $\sigma$.
Then $P_{\tau}$ is a smooth point of $X_{\sigma}$ if and only if
 $P_{\tau}$ is a smooth point of $X_{\tau}$, i.e., if and only if
  $\tau$  is generated by a part of a basis of $N$
($N$ is the $\Z$ dual of the character group of the torus).
\end{lemma}
\begin{proof} We have that $X_{\tau}$ is a principal open subset of $X_{\sigma}$.
  Hence $X_\sigma$ is non-singular at $P_\tau$ if and only if
  $X_\tau$ is non-singular at $P_\tau$.  By Lemma~\ref{l-1-tau},
  $X_\tau$ is non-singular at $P_\tau$ if and only if $X_\tau$ is a
  non-singular variety; but by \S 2.1 of \cite{F}, this is true if and only if $\tau$ is generated by a part of a basis of $N$.
\end{proof}

\begin{theorem}\label{sing} Let $\tau =\sigma_{i,j}$. We have

\begin{enumerate}\item $P_\tau\in $Sing$\, X_\sigma$.
\item We have an identification of $TC_{P_{\tau}}X_\sigma$ with
$Z\times (K)^{\# J({\L})-3}, Z$ being as in Theorem \ref{quadric};
in particular, $TC_{P_{\tau}}X_\sigma$ is a toric variety. \item
The singularity at $P_\tau$ is of the same type as that at the
vertex of the cone over the quadric surface $x_1x_4-x_2x_3=0$ in
$\mathbb{P}^3$. In particular, mult$_{P_{\tau}}X_\sigma=2$.
\end{enumerate}
\end{theorem}

\begin{proof}
Assertion (1) follows from Lemmas \ref{l-1-tau}, \ref{basis-sing}
and Theorem \ref{quadric}.

 \ni $X_\tau$ being open in $X_\sigma$, we may identify
$TC_{P_{\tau}}X_\sigma$ with $TC_{P_{\tau}}X_\tau$; assertion (2)
follows from this in view of Theorem \ref{quadric} (and \S
\ref{hilbert}).

\ni Assertion (3) is immediate from (2).
\end{proof}

Next, we will show that the faces containing some $\sigma_{ij}$
are the only singular faces. We first prove some preparatory
Lemmas.

\begin{lemma}\label{3} Let $A\not= {\hat{0}}$.
 If $e_A-e_C$ is in $W$ (the set of generators of $\sigma$ as
 described in Proposition \ref{prop-cpc}),
 then it is in $W(\sigma_{ij})$ (cf. Definition \ref{gen})
 for some $(i,j)$ where $1\leq i \leq n-d-1$, $1\leq j\leq d-1$.
 \end{lemma}
\begin{proof} If $A$ is equal to some $\mu_{ij}$,
then $C$ must be one of the two covers of
$\mu_{ij}(=(1,\cdots,j,i+j+1,\cdots,i+d))$ in $J({\mathcal{L}})$,
and we are done by Theorem~\ref{main}. Thus, we will assume that
$A\not=\mu_{ij}$. Hence $A$ is a join irreducible of one of the
following two forms.

\noindent\textbf{Case 1:} $A=(1,\ldots , k , n-d+k+1,\ldots , n)$ for
some $k$.

Then $\mu_{n-d-1,k}=(1,\ldots , k,n-d+k,\ldots , n-1)$; and
$\left(A,\mu_{n-d-1,k}\right)$ is a cover in $J\left({\L}\right)$.
Also, $A$ only has one cover in $J\left({\L}\right)$, which must
be $C$, thus $e_A-e_C$ is an element of
$W\left(\sigma_{n-d-1,k}\right)$ as shown in Cases 1 and 2 of
Theorem~\ref{main}.

\noindent\textbf{Case 2:} $A=(k+1,\ldots , k+d)$, $1\leq k\leq
n-d-1$ (note that $k<n-d$, since $C>A$ because $e_A-e_C\in W$).

Then we have $\mu_{k,1}=(1,k+2,\ldots , k+d)$; and $(A,\mu_{k,1})$
is a cover in $J\left({\L}\right)$.  Also, we must have
$C=(k+2,\ldots , k+d+1)$, and $e_A-e_C$ is an element of
$W\left(\sigma_{k,1}\right)$ by Case 1 of Theorem~\ref{main}.
\end{proof}

We now return to the case of a G-H toric variety.

\begin{theorem}\label{no-loops-allowed} Let $\tau$ be a face such that
 $D_{\tau}$ is not contained in any $\Lij$,
 $1\leq i\leq n-d-1$, $1\leq j\leq d-1$. Then the associated face $\tau$ is nonsingular (i.e., if any face $\tau$ does not contain any one $\sigma_{ij}$, then $\tau$ is nonsingular).
\end{theorem}
\begin{proof}  By Lemma \ref{basis-sing}, for $\tau$
to be nonsingular, it must be generated by part of a basis for
$N$.  Since $\tau$ is generated by a subset $W(\tau)$ of $W$, for
$\tau$ to be singular its generators would have to be linearly
dependent. (Generally this is not enough to prove a face is
singular or nonsingular, but since all generators in $W$ have
coefficients equal to $\pm 1$, any linearly independent set will
serve as part of a basis for $N$.) Suppose $\tau$ is singular;
then there is some subset of the elements of $W\left(\tau\right)$
equal to $\{e_1-e_2,\ldots\}$ such that $\sum
a_{ij}\left(e_i-e_j\right)=0$, with coefficients $a_{ij}$ nonzero
for at least one $(i,j)$.

Recall that the elements of $W$ can be represented as all the line
segments in the lattice $J\left({\L}\right)$, with the exception
of $e_{\hat{1}}$ (see diagram in Remark~\ref{diagram-remark}).
Therefore, the linearly dependent generators of $\tau$ must
represent a ``loop''  of line segments in $J\left({\L}\right)$.
This loop will have at least one bottom corner, left corner, top
corner, and right corner.

Choose some particular $\Lij$.  From the previous section, we have
that
$W\left(\sigma_{ij}\right)=\{e_{\mu_{ij}}-e_A,\,e_{\mu_{ij}}-e_B,\,e_A-e_C,\,e_B-e_C\}$.
These four generators are represented by the four sides of a
diamond in $J\left({\L}\right)$.   Thus, by hypothesis, the
generators of $\tau$ represent a loop in $J\left({\L}\right)$ that
does not traverse all four sides of the diamond representing all
four generators of $\sigma_{ij}$.

Since by hypothesis, $D_{\tau}$ is not contained in any $\Lij,1\le
i\le n-d-1,1\le j\le d-1$, there must be at least one element of
$D_{\tau}$ in the interval $[\mu_{ij},\lambda_{ij}]$, say
$\alpha\in [\mu_{ij},\lambda_{ij}]$.  We have,
$\alpha\geq\mu_{ij}$ and $\alpha\not\geq C$ for $C$ as defined in
the proof of Theorem~\ref{main}.  Based on how $\alpha$ compares
to both $A$ and $B$, we can eliminate certain elements of $W$ from
$W\left(\tau\right)$.  There are four possibilities; we list all
four, as well as the corresponding generators in $W(\sigma_{ij})$
which are not in $W(\tau)$, i.e., those generators $v$ in
$W(\sigma_{ij})$ such that $f_{I\left(\alpha\right)}(v)\not= 0$:

\begin{eqnarray*}
\alpha\not\geq A,\,\alpha\not\geq B &\Rightarrow & e_{\mu_{ij}}-e_A,\,e_{\mu_{ij}}-e_B\not\in W(\tau)\\
\alpha\geq A,\,\alpha\not\geq B & \Rightarrow &
e_A-e_C,\,e_{\mu_{ij}}-e_B\not\in W(\tau)\\
\alpha\not\geq A,\,\alpha\geq B & \Rightarrow & e_{\mu_{ij}}-e_A,\,e_B-e_C\not\in W(\tau)\\
\alpha\geq A,\,\alpha\geq B & \Rightarrow &
e_A-e_C,\,e_B-e_C\not\in W(\tau)
\end{eqnarray*}

Therefore, it is impossible to have $\{e_{\mu_{ij}}-e_A,e_A-e_C\}$
or $\{e_{\mu_{ij}}-e_B,e_B-e_C\}$ contained in $W(\tau)$.  This is
true for any $(i,j)$; therefore (in view of Lemma~\ref{3}) our
``loop'' in $J\left({\L}\right)$ that represented the generators
of $\tau$ cannot have a left corner or a right corner. Therefore
it is really not possible to have a loop at all, and so the
generators of $\tau$ are linearly independent, and the result
follows.
\end{proof}

\begin{cor}\label{nonsing} The G-H toric variety $X_{d,n}$ is smooth along the orbit $O_{\tau}$ if and only if the face $\tau$ does not contain any $\sigma_{ij}$.
\end{cor}
Combining the above Theorem with Theorem \ref{sing} and Lemma
\ref{dim3}, we obtain our first main Theorem:
\begin{theorem}\label{mainth} Let ${\L}=I_{d,n}$.
 Then
\begin{enumerate}
 \item $\displaystyle Sing\,
X\left(\L\right)=\bigcup_{\sigma_{i,j}}
\overline{O}_{\sigma_{i,j}}$, the union being taken over all
$\sigma_{i,j}$'s as in Theorem \ref{main}.
 \item Sing\,$X\left(\L\right)$ is
pure of codimension $3$ in $X\left(\L\right)$; further, the
generic singularities are of cone type (more precisely, the
singularity type is the same as that at the vertex of the cone over
the quadric surface $x_1x_4-x_2x_3=0$ in $\mathbb{P}^3$ ). \item
For $\tau=\sigma_{i,j}$, $TC_{P_{\tau}}X({\L})$ is a toric
variety; further, mult$_{P_{\tau}}X\left(\L\right)=2$.
\end{enumerate}
\end{theorem}

\begin{remark} Thus Theorem~\ref{mainth} proves the
conjecture of \cite{g-l} using just the combinatorics of the cone
associated to the toric variety $X_{d,n}$ (for a statement of the
conjecture of \cite{g-l}, see Remark \ref{conj11}). Further, it
gives a description of Sing$\,X_{d,n}$ purely in terms of the
faces of the cone associated to $X_{d,n}$.
\end{remark}


\section{Multiplicities of Singular Faces of $X_{2,n}$}\label{d=2-sect}

In this section, we take ${\mathcal{L}}=I_{2,n}$, determine the
multiplicity of $X_{2,n} (=X\left(I_{2,n}\right))$ at $P_\tau$ for
certain of the singular faces of $X_{2,n}$, and deduce a product
formula. For $I_{d,n}$, we defined $\Lij$ and the corresponding
face $\sigma_{i,j}$ for $1\leq j\leq d-1$, $1\leq i\leq n-d-1$;
thus for $I_{2,n}$, we need to consider only $\mathcal{L}_{i,1}$
for $1\leq i\leq n-3$.

For example, below is the poset of join irreducibles for
$I_{2,6}$.  We write $\sigma_{i,1}$ inside each diamond because
the four segments surrounding it represent the four generators of
the face.

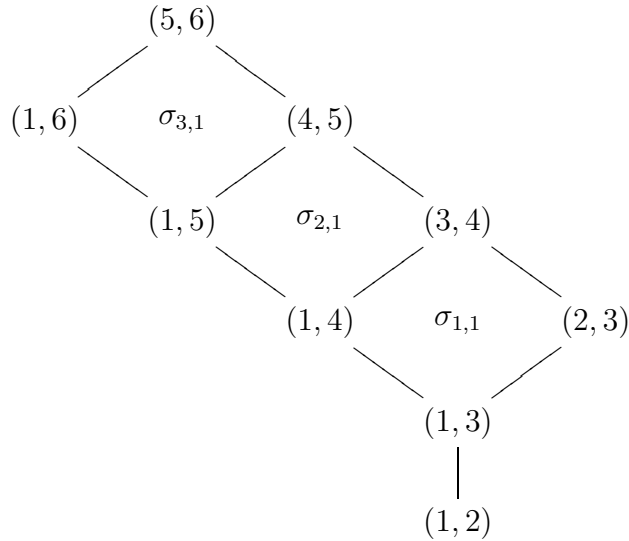
\[ \xymatrix@-5pt{
 & & (5,6) & & & \\
 &(1,6)\ar@{-}[ur] & \sigma_{3,1} & (4,5)\ar@{-}[ul] & & \\
 & & (1,5)\ar@{-}[ul]\ar@{-}[ur] & \sigma_{2,1} & (3,4)\ar@{-}[ul] & \\
 & & & (1,4)\ar@{-}[ul]\ar@{-}[ur]& \sigma_{1,1} & (2,3)\ar@{-}[ul] \\
 & & & & (1,3)\ar@{-}[ul]\ar@{-}[ur] & \\
 & & & & (1,2)\ar@{-}[u] &
 } \]

To go from the join irreducibles of $I_{2,6}$ to $I_{2,7}$, we
just add $(1,7)$ and $(6,7)$ to the poset above, forming
$\sigma_{4,1}$.  We will see that this makes the calculation of
the multiplicities of singular faces of $I_{2,n}$ much easier.

In the sequel, we shall denote the set of join irreducibles of
$I_{2,n}$ by $J_{2,n}$; also, as in the previous sections,
$\sigma$ will denote the cone corresponding to $X_{2,n}$.

\subsection{Mult$_{P_\sigma}\, X_{2,n}$}\label{mult-sigma}

Let $\sigma$ be the CPC associated to $I_{d,n}$. Now $X_{d,n}$
being of cone type (i.e., the vanishing ideal is homogeneous), we
have a canonical identification of $T_{P_{\sigma}}X_{d,n}$ (the
tangent cone to $X_{d,n}$ at $P_\sigma$) with $X_{d,n}$. Hence by
Theorem \ref{DIM'},(2), we have that mult$_{P_\sigma}\, X_{d,n}$
equals the number of maximal chains in $I_{d,n}$. So we begin by
counting the number of maximal chains in $I_{2,n}$.

As we move through a chain from $(1,2)$, at any point $(i,j)$ we
have at most two possibilities for the next point: $(i+1,j)$ or
$(i,j+1)$.  For each cover in our chain, we assign a value: for a
cover of type $\left(\left(i,j+1\right),\left(i,j\right)\right)$
assign $+1$; for a cover of type
$\left(\left(i+1,j\right),\left(i,j\right)\right)$ assign $-1$.

A maximal chain $C$ in $I_{2,n}$ contains $2n-3$ lattice points,
and thus every chain can be uniquely represented by a
$(2n-4)$-tuple of $1$'s and $-1$'s; let us denote this
$(2n-4)$-tuple by $n_C=\left<a_1,\ldots , a_{2n-4}\right>$.

What rules do we have on $n_C$?  First it is clear that $1$ and
$-1$ occur precisely $n-2$ times.  Secondly, we can see that
$a_1=+1$, and for any $1\leq k\leq 2n-4$, if $\{a_1,\ldots ,a_k\}$
contains more $-1$'s than $+1$'s, then we have arrived at a point
$(i,j)$ with $i>j$, which is not a lattice point.  Thus, we must
have $a_1+\ldots +a_k\geq 0$ for every $1\leq k\leq 2n-4$.

\begin{theorem}[cf. Corollary 6.2.3 in \cite{stanley}]\label{cat} The Catalan number
\[ Cat_n=\frac{1}{n+1}{2n\choose n},\ (n\geq 0)\]
counts the number of sequences $a_1,\ldots , a_{2n}$ of $1$'s and
$-1$'s with
\[a_1+\ldots +a_k\geq 0,\  (k=1,2,\ldots , 2n)\]
and $a_1+\ldots +a_{2n}=0$.
\end{theorem}

\begin{cor}\label{cor-8-3} The multiplicity of $X_{2,n}$  at $P_\sigma$ is equal to the Catalan number
\[Cat_{n-2}=\frac{1}{n-1}{2n-4\choose n-2}.\]
\end{cor}

\subsection{Mult$_{P_\tau}\,X_{2,n}$}\label{mult-tau}

In this section, we determine mult$_{P_\tau}\,X_{2,n}$, for $\tau$
being of block type (see Definition \ref{block} below). Let $\tau$
be a face of $\sigma$, such that the associated (embedded
sublattice) $D_\tau$ is of the form:
\[D_\tau = [(1,2),(i,i+1)]\cup [(i+k+2,i+k+3),(n-1,n)]=I_1\cup I_2,\mathrm{\ say},\]
where $I_1=[(1,2),(i,i+1)],$ $I_2=[(i+k+2,i+k+3),(n-1,n)],$ $1\le i\le
n-3,$ $0\le k\le n-i-3$.

We shall now determine $W(\tau)$ (cf. \S \ref{gen}). Let $A_\tau$
denote the interval

\ni $[(1,i+2),(i+k+2,i+k+3)]$ in $J_{2,n}$:
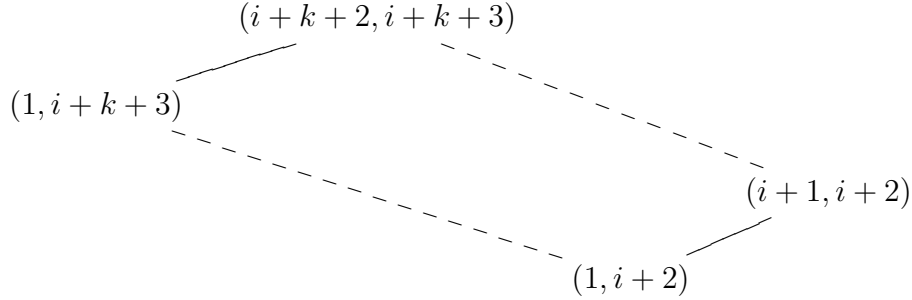
\[ \xymatrix@-10pt{
 & (i+k+2,i+k+3) & &  \\
 (1,i+k+3)\ar@{-}[ur] & & &  \\
 & & & (i+1,i+2)\ar@{--}[uull] \\
 & & (1,i+2)\ar@{--}[uull]\ar@{-}[ur] & \\
 } \]

\begin{lem}\label{asin}
With $\tau$ as above, we have
$W(\tau)=\{e_{y'}-e_y\,|\,(y,y')\mathrm{\ is\ a\ cover\ in\
}A_\tau\}$.
\end{lem}
\begin{proof} Clearly $e_{(n-1,n)}$ (the element in
$W(\sigma)$ corresponding to the unique maximal element $(n-1,n)$
in $J_{2,n}$) is not in $W(\tau)$ (since $(n-1,n)\in D_\tau$). Let
us denote
\[\theta=(i+k+2,i+k+3),\delta=(i.i+1)\]

\vs.1cm\ni\textbf{Claim 1:} For a cover $(y,y')$ in $A_\tau$,
$f_{I_{\alpha}}(e_{y'}-e_y)=0, \forall\alpha\in D_\tau$.

\ni The claim follows in view of the facts that for a cover
$(y,y')$ in $A_\tau$, we have,

(i) $y,y'\in I_\theta$, and hence $y,y'\in I_\alpha,
\forall\alpha\in I_2$.

(ii) $y,y'\not\in I_\delta$, and hence $y,y'\not\in I_\alpha,
\forall\alpha\in I_1$.

\vs.1cm\ni\textbf{Claim 2:} For a cover $(y,y')$ in $J_{2,n}$ not
contained in $A_\tau$, there exists an $\alpha\in D_\tau$ such
that $f_{I_{\alpha}}(e_{y'}-e_y)\ne 0$.

\ni Note that a cover in $J_{2,n}$ is one of the following three
types:

\textbf{Type I:} $\left((1,j), (1,j-1)\right), 3\le j\le n$

\textbf{Type II:} $\left((j-1,j),(j-2,j-1)\right), 4\le j\le n$

\textbf{Type III:} $\left((j-1,j), (1,j)\right), 3\le j\le n$

\ni Let now $(y,y')$ be a cover not contained in $A_\tau$.

If $(y,y')$ is of Type I, then $(y,y')=\left((1,j),
(1,j-1)\right)$, where either $j\le i+2$ or $j\ge i+k+4$. Letting
 $$\alpha=\begin{cases}(1,j-1), & {\mathrm{\ if\ } j\le i+2}\\
 (j-2,j-1), & \mathrm{\ if\ } j\ge i+k+4
  \end{cases}$$ we have, $\alpha\in D_\tau$ and
$f_{I_{\alpha}}(e_{y'}-e_y)\ne 0$.

If $(y,y')$ is of Type II, then
$(y,y')=\left((j-1,j),(j-2,j-1)\right)$, where either $j\le i+2$
or $j\ge i+k+4$. Letting
 $\alpha=(j-2,j-1)$, we have, $\alpha\in D_\tau$ and
$f_{I_{\alpha}}(e_{y'}-e_y)\ne 0$.

If $(y,y')$ is of Type III, then $(y,y')=\left((j-1,j),
(1,j)\right)$, where either $j\le i+1$ or $j\ge i+k+4$. Letting
 $$\alpha=\begin{cases}(1,j), & {\mathrm{\ if\ } j\le i+1}\\
 (j-2,j), & \mathrm{\ if\ } j\ge i+k+4
  \end{cases}$$ we have, $\alpha\in D_\tau$ and
$f_{I_{\alpha}}(e_{y'}-e_y)\ne 0$.

The required result follows from claims 1 and 2.
\end{proof}
\begin{cor}\label{btype}
With $\tau$ as in Lemma \ref{asin}, we have
\[\tau = \sigma _{i,1}\cup\sigma_{i+1,1}\cup\ldots\cup\sigma_{i+k,1}.\]
\end{cor}
\begin{definition}\label{block} We define a face $\tau$ as in
 Lemma \ref{asin}
 as a \emph{$J$-block} (namely, $\tau$ is an union of consecutive
 $\sigma_{i,1}$'s).
 \end{definition}

 \begin{rem} Note that in general an union of faces need not be
  a face.
\end{rem}
\subsection{The Hibi variety $Z_{2,r}$}\label{z2r} For an integer $r\ge 3$,
let ${\widetilde{I_{2,r}}}$ denote the distributive lattice
$I_{2,r}\,\setminus\,\{(1,2),(r-1,r)\}$. We define $Z_{2,r}$ to be
the Hibi variety associated to ${\widetilde{I_{2,r}}}$. Note (cf.
Proposition \ref{prop-cpc}) that the cone associated to $Z_{2,r}$
has a set of generators consisting of $\{e_{y'}-e_y\}$, $(y,y')$
being a cover in the sublattice $[(1,3),(r-1,r)]$ of $J_{2,r}$
(the set of join irreducibles of $I_{2,r}$). In view of Theorem
\ref{DIM'}, (2), we have
\[\mathrm{mult}_{\mathbf{0}}Z_{2,r}=\mathrm{mult}_{\mathbf{0}}X_{2,r}
(=Cat_{r-2})\] (Here, $\mathbf{0}$ denotes the origin.)
\begin{theorem}\label{z2r'} Let $\tau$ be a face of $\sigma$ which
 is a ``$J$-block'' of $k+1$ consecutive $\sigma_{i,1}$'s
 (as in Definition \ref{block}).
   We have an identification of $X_\tau$ (the open affine piece
of $X_\sigma$ corresponding to $\tau$) with $Z_{2,k+4}\times
(K^*)^m$ where $m=\mathrm{codim}_\sigma\tau(=2(n-k)-6)$
\end{theorem}
\begin{proof}  In view of \ref{orbit-decomp} and
Proposition \ref{orbitdim}, we have

\ni
$\mathrm{codim}_\sigma\tau=dim\,X(D_\tau)=\#\{\mathrm{elements\
in\ a\ maximal\ chain\ in\ }D_\tau\}$; from this it is clear that
$\mathrm{codim}_\sigma\tau $ equals $2(n-k)-6$. Next, in view of
Lemma \ref{asin} and \S \ref{z2r}, we obtain an identification of
$X_\tau$ with $Z_{2,k+4}\times (K^*)^m$ ($m$ being as in the
Theorem).
\end{proof}

\begin{theorem}\label{cat-mult} Let $\tau$ be as
in Theorem \ref{z2r'}.
  \begin{enumerate}
\item We have an identification of $TC_{P_{\tau}}X_\sigma$ with
$Z_{2,k+4}\times (K)^m$ where

\ni $m=\mathrm{codim}_\sigma\tau(=2(n-k)-6)$;
  in particular, $TC_{P_{\tau}}X_\sigma$ is a toric variety.
  \item  mult$_{P_{\tau}}X_{2,n}=Cat_{k+2}
  (=\frac{1}{k+3}{2k+4\choose k+2})$.
\end{enumerate}
\end{theorem}
\begin{proof} $X_\tau$ being open in $X_\sigma$, we may identify
$TC_{P_{\tau}}X_\sigma$ with $TC_{P_{\tau}}X_\tau$; assertion (1)
follows from this in view of Theorem \ref{z2r'} (and \S
\ref{hilbert}).

\ni Assertion (2) follows from (1) and Corollary~\ref{cor-8-3}.
\end{proof}

\subsection{A Product Formula}\label{product-formula-subsection}

In this subsection, we give a product formula for

\ni $mult_{P_\tau} \,X_{2,n}$ where $\tau$ is a union of pairwise
non-intersecting and non-consecutive $J$-blocks (see Remark
\ref{see} below).

Let $\tau$ be a face of $\sigma$, such that the associated
(embedded sublattice) $D_\tau$ is of the form:
\[D_{\tau}=[(1,2),(i_1,i_1+1)]\cup [(i_1+k_1+2,i_1+k_1+3),(i_2,i_2+1)]\]
\[\cup [(i_2+k_2+2,i_2+k_2+3),(n-1,n)]=J_1\cup J_2\cup J_3,\,
\mathrm{say}\] where $i_1+k_1+1< i_2$, and
$J_1=[(1,2),(i_1,i_1+1)],J_2=[(i_1+k_1+2,i_1+k_1+3),(i_2,i_2+1)]$

\ni $J_3=[(i_2+k_2+2,i_2+k_2+3),(n-1,n)]$.

Consider the following sublattices in $J_{2,n}$ (the set of
join-irreducibles in $I_{2,n}$):
\[A=[(1,i_1+2),(i_1+k_1+2,i_1+k_1+3)], B=[(1,i_2+2),(i_2+k_2+2,i_2+k_2+3)] . \]
\begin{lem}\label{asin'}
With $\tau$ as above, we have
$W(\tau)=\{e_{y'}-e_y\,|\,(y,y')\mathrm{\ is\ a\ cover\ in\ }A\cup
B\}$.
\end{lem}
\begin{proof} We proceed as in the proof of Lemma \ref{asin}. As in that proof,
we have $e_{(n-1,n)}$ is not in $W(\tau)$ (since $(n-1,n)\in
D_\tau$). Let us denote
\[\theta_1=(i_1+k_1+2,i_1+k_1+3),\theta_2=(i_2+k_2+2,i_2+k_2+3),
\delta_1=(i_1,i_1+1),\delta_2=(i_2,i_2+1) . \] For any cover $(y,y')$
in $A\cup B$, we clearly have $y,y'\in I_{\theta_2}$ and hence
$y,y'\in I_\alpha,\forall\alpha\in J_3$; also, $y,y'\not\in
I_{\delta_1}$ and hence  $y,y'\not\in I_\alpha,\forall\alpha\in
J_1$. Thus we obtain that
$$f_{I_{\alpha}}(e_{y'}-e_y)=0,\forall\alpha\in J_1\cup
J_3\leqno{(*)}$$

Next, if $(y,y')$ is a cover in $A$, then $y,y'\in I_{\theta_1}$
and hence $y,y'\in I_\alpha,\forall\alpha\in J_2$; also, if
$(y,y')$ is a cover in $B$, then $y,y'\not\in I_{\delta_2}$, and
hence $y,y'\not\in I_\alpha,\forall\alpha\in J_2$ (note that
$\theta_1$ (resp. $\delta_2$) is the smallest (resp. largest)
element in $J_2$). Thus we obtain that
$$f_{I_{\alpha}}(e_{y'}-e_y)=0,\forall\alpha\in J_2\leqno{(**)}$$

Now $(*)$ and $(**)$ imply the inclusion ``$\supseteq$". We shall
prove the inclusion ``$\subseteq$" by showing that if a cover
$(y,y')$ is not contained in $A\cup B$, then there exists an
$\alpha\in D_\tau$ such that $f_{I_{\alpha}}(e_{y'}-e_y)\ne 0$.
This proof again runs on similar lines as that of Lemma
\ref{asin}. Let then $(y,y')$ be a cover in $J_{2,n}$ not
contained in $A\cup B$. It would be convenient to introduce the
following sublattices in $J_{2,n}$:
\begin{eqnarray*}P & = & [(1,2),(i_1+1,i_1+2)]\\
Q & = & [(1,i_1+k_1+3),(i_2+1,i_2+2)]\\
R & = & [(1,i_2+k_2+3),(n-1,n)].
\end{eqnarray*}

We distinguish the following cases:

\ni\textbf{Case I:} Let $(y,y')$ be of type I (cf. proof of Lemma
\ref{asin}), say, $\left((1,j), (1,j-1)\right)$.

(i) If $(y,y')$ is contained in $P$, then $ j\le i_1+2$. We let
$\alpha=(1,j-1)$. Note that $\alpha\in J_1$ and
$f_{I_{\alpha}}(e_{y'}-e_y)\ne 0$.

(ii) If $(y,y')$ is contained in $Q$ (resp. $R$), then $
i_1+k_1+4\le j\le i_2+2$ (resp.$i_2+k_2+4\le j\le n$). We let
$\alpha=(j-2,j-1)$. Note that $\alpha\in J_2$ (resp. $J_3$) and
$f_{I_{\alpha}}(e_{y'}-e_y)\ne 0$.

\ni\textbf{Case II:} Let $(y,y')$ be of type II, say,
$\left((j-1,j), (j-2,j-1)\right)$.

Then, $ 3\le j\le i_1+2$, $ i_1+k_1+4\le j\le i_2+2$ or $
i_2+k_2+4\le j\le n$ accordingly as $(y,y')$ is contained in $P$,
$Q$ or $R$. We let $\alpha=(j-2,j-1)$. Note that $\alpha\in J_1$,
$J_2$ or $J_3$, accordingly as $(y,y')$ is contained in $P$, $Q$ or
$R$; and $f_{I_{\alpha}}(e_{y'}-e_y)\ne 0$.

\ni\textbf{Case III:} Let $(y,y')$ be of type III, say,
$\left((j-1,j), (1,j)\right)$.

(i) If $(y,y')$ is contained in $P$, then $ j\le i_1+1$. We let
$\alpha=(1,j)$. Note that $\alpha\in J_1$ and
$f_{I_{\alpha}}(e_{y'}-e_y)\ne 0$.

(ii) If $(y,y')$ is contained in $Q$ (resp. $R$), then $
i_1+k_1+4\le j\le i_2+1$ (resp. $ i_2+k_2+4\le j\le n$). We let
$\alpha=(j-2,j)$. Note that $\alpha\in J_2$ or $J_3$, accordingly as
$(y,y')$ is contained in $Q$ or $R$; and
$f_{I_{\alpha}}(e_{y'}-e_y)\ne 0$.
\end{proof}

As an immediate consequence of Lemma \ref{asin'} and Corollary
\ref{btype}, we have

\begin{cor}\label{non}
Let $\tau$ be as in Lemma \ref{asin'}. Then
$\tau=\tau_1\cup\tau_2$, where
\[\tau_1=\sigma_{i_1,1}\cup\ldots\cup \sigma_{i_1+k_1,1},\]
\[\tau_2=\sigma_{i_2,1}\cup\ldots\cup \sigma_{i_2+k_2,1};\]
\end{cor}
\begin{rem}\label{see}
We refer to a pair $(\tau_1,\tau_2)$ of faces as in Corollary
\ref{non} as \emph{non-intersecting $J$-blocks}.
\end{rem}

\begin{theorem}\label{z2r''} Let $\tau=\tau_1\cup\tau_2$, where $\tau_1$ and $\tau_2$ are two
non-intersecting (and non-consecutive) $J$-blocks (cf. Corollary
\ref{non}). We have an identification of $X_\tau$ (the open affine
piece of $X_\sigma$ corresponding to $\tau$) with
$Z_{2,k_1+4}\times Z_{2,k_2+4}\times (K^*)^m$

\ni where $m=\mathrm{codim}_\sigma\tau(=2(n-k_1-k_2)-9)$
\end{theorem}

Proof is similar to that of Theorem \ref{cat-mult} (using Lemma
\ref{asin'}).

\ni As an immediate consequence, we have
\begin{theorem}\label{cat-mult''} Let $\tau=\tau_1\cup\tau_2$, where $\tau_1$ and $\tau_2$ are two
non-intersecting (and non-consecutive) $J$-blocks (cf. Corollary
\ref{non})
  Then
  \begin{enumerate}

\item We have an identification of $TC_{P_{\tau}}X_\sigma$ with

\ni$Z_{2,k_1+4}\times Z_{2,k_2+4}\times (K^*)^m$ where
$m=\mathrm{codim}_\sigma\tau(=2(n-k_1-k_2)-9)$; in particular,
$TC_{P_{\tau}}X_\sigma$ is a toric variety.
  \item  $mult_{P_{\tau}}X_{2,n}=
  \left(mult_{P_{\tau_1}}X_{2,n}\right)
  \cdot\left(mult_{P_{\tau_2}}X_{2,n}\right)$.
\end{enumerate}
\end{theorem}

Proof is similar to that of Theorem \ref{cat-mult} (using Theorem
\ref{z2r''}).


\begin{remark} It is clear that we can extend this multiplicative property to $\tau=\tau_1\cup\ldots\cup\tau_s$, a union of $s$ pairwise non-intersecting, non-consecutive $J$-blocks.
\end{remark}

\section{A Multiplicity Formula for $X_{\MakeLowercase{d,n}}$}\label{last-sect}

In this section, we give a formula for $mult_{P_\sigma}X_{d,n}$.
 By Theorem \ref{DIM'}, (2), we have that $mult_{P_\sigma}X_{d,n}$ equals the
number of maximal chains in $I_{d,n}$.  We give below an explicit
formula for the number of maximal chains in $I_{d,n}$.  Notice
that the number of chains in $I_{d,n}$ from $(1,2,\ldots , d)$ to
$(n-d+1,\ldots , n)$ is the same as the number of chains from
$(0,0,\ldots 0)$ to $(n-d,n-d,\ldots n-d)$ such that for any
$(i_1,\ldots i_d)$ in the chain, we have $i_1\geq i_2\geq \ldots
i_d\geq 0$.  Now, set
\[\mu=(\mu_1,\mu_2,\ldots ,\mu_d )=(n-d,n-d,\ldots n-d).\eqno{(\ast)}\]

For any $\lambda\vdash m$, let $f^{\lambda}=K_{\lambda, 1^{m}}$,
i.e., the number of standard Young tableau of shape $\lambda$ (cf.
\cite{stanley}).

\begin{prop}[cf. Proposition 7.10.3, \cite{stanley}] Let $\lambda$ be a partition of $m$.  Then the number $f^\lambda$ counts the lattice paths $0=v_0,v_1,\ldots v_{m}$ in $\R^l$ (where $l=l\left(\lambda\right)$) from the origin $v_0$ to $v_{m}=(\lambda_1,\lambda_2,\ldots \lambda_l)$, with each step a coordinate vector; and staying within the region (or cone) $x_1\geq x_2\geq \ldots \geq x_l\geq 0$.
\end{prop}

Thus, for $\mu$ as described in $(\ast)$ above, we have that the
number of maximal chains in $I_{d,n}$ is equal to $f^{\mu}$.

Corollary 7.21.5 of \cite{stanley} gives an explicit description
of $f^{\lambda}$.

\begin{prop} Let $\lambda\vdash m$.  Then
\[ f^{\lambda}=\frac{m!}{\prod _{u\in\lambda}h\left(u\right)}.\]
\end{prop}

The statement above refers to $u\in\lambda$ as a box in the Young
tableau of $\lambda$, and $h\left(u\right)$ being the ``hook
length'' of $u$.  The hook length is easily defined as the number
of boxes to the right, and below, of $u$, including $u$ once.

Let us take, for example, $I_{3,6}$.  Then $\mu =(3,3,3)$, and the
Young tableau of shape $\mu$ with hook lengths given in their
corresponding boxes is
\begin{center}
\setlength{\unitlength}{.3 in}
\begin{picture}(3, 3)
 \put(0,2){\framebox(1,1){5}}
 \put(0,1){\framebox(1,1){4}}
 \put(0,0){\framebox(1,1){3}}
 \put(1,2){\framebox(1,1){4}}
 \put(1,1){\framebox(1,1){3}}
 \put(1,0){\framebox(1,1){2}}
 \put(2,2){\framebox(1,1){3}}
 \put(2,1){\framebox(1,1){2}}
 \put(2,0){\framebox(1,1){1}}
\end{picture}
\end{center}
Therefore
\[f^{\mu}=\frac{9!}{5\cdot 4^2\cdot 3^3\cdot 2^2\cdot 1}=42.\]

In fact, in the scenario of $I_{d,n}$ our derived partition $\mu$
(given by ($\ast$)) will always be a rectangle; and we can deduce
a formula for $f^{\mu}$ which does not require the Young tableau.
The top left box of $\mu$ will always have hook length
$(n-d)+d-1=n-1$. Then, the box directly below it, and the box
directly to the right of it will have length $n-2$.  For any box
of $\mu$, the box below and the box to the right will have hook
length one less than that of the box with which we started.

Since the posets $I_{d,n}$, and $I_{n-d,n}$ are isomorphic, we can
assume that $d\leq n-d$.  Then we have $\displaystyle
\prod_{u\in\mu}h\left(u\right)= $
\[(n-1)(n-2)^2\cdots (n-d)^d(n-d-1)^d\cdots (d)^d(d-1)^{d-1}\cdots (2)^2(1).\]
Thus we arrive at the following.

\begin{theorem} The multiplicity of $X_{d,n}$ at $P_{\sigma}$ is equal to
\[\frac{\left(d\left(n-d\right)\right)!}{(n-1)(n-2)^2\cdots (n-d)^d(n-d-1)^d\cdots (d)^d(d-1)^{d-1}\cdots (2)^2(1)}.\]
\end{theorem}
\section{Conjectures}\label{conjectures} In this section, we give two conjectures on
the multiplicity at a singular point and one conjecture on
Sing$\,X({\L}),{\L}$ being the Bruhat poset of Schubert varieties
in a minuscule $G/P$.

Guided by the phenomenon for $X_{2,n}$ as given by
 Theorem~\ref{cat-mult''}, we make the following conjecture.\\

\noindent\textbf{Conjecture 1.} There exists a multiplicative
formula in $X_{d,n}$ similar to Theorem \ref{cat-mult''} for a
face which has the form as an union of several disjoint
(and non-consecutive) faces.\\

The generating set $W(\tau)$ of a face $\tau$ consists of
$\{e_{y'}-e_y\}$, for certain covers $(y,y')$ in $J(\mathcal{L})$
(assuming that ${\hat{1}}\in D_\tau$ so that $e_{{\hat{1}}}$ is
not in $W(\tau)$). Thus $W(\tau)$ determines a sub set
$H(\tau):=\dot\bigcup\,H(\tau)_{i}$ of $J(\mathcal{L})$, such that
$W(\tau)$ consists of all the covers in the $H(\tau)_{i}$'s. Thus,
$H(\tau)$ for $\tau = \sigma_{ij}$ would be the diamond given in
Remark~\ref{sigma-ij-diamond}. In \S
\ref{product-formula-subsection}, if $\tau = \tau_1\cup \tau_2$
for $\tau_1,\tau_2$ a pair of non-consecutive, non-intersecting
$J$-blocks, $H(\tau ) = H(\tau_1 )\dot\cup H(\tau_2 )$.

Theorem~\ref{cat-mult} implies that $mult_{P_{\tau_1}}\,X_{2,n}=
mult_{P_{\tau_2}}\,X_{2,n}$ if both $\tau_1$ and $\tau_2$ are
$J$-blocks of the same length; in particular, $H(\tau_1)$ and
$H(\tau_2)$ are isomorphic.  Guided by this phenomenon we make the following conjecture.\\

\noindent\textbf{Conjecture 2.} For a face $\tau$ of any Hibi
toric variety $X\left(\L\right)$, $mult_{P_\tau}\,
X\left(\L\right)$ is determined by the poset $H(\tau)$.  By this
we mean that if $\tau$, $\tau '$ are such that $H(\tau)$,
$H(\tau')$ are isomorphic posets, then the multiplicities of
$X\left(\L\right)$ at the points
$P_{\tau}$, $P_{\tau '}$ are equal.\\

\begin{remark}\label{conj11} Toward generalizing
Theorem \ref{mainth} to other Hibi varieties, we will first
explain how the lattice points $\mu_{ij}$, $\lambda_{ij}$ were
chosen. Let $\theta$, $\delta$ be two incomparable, join and meet
irreducibles  in $I_{d,n}$; say $\theta = (1,\ldots j,
n+j+1-d,\ldots n)$ and $\delta = (i+1,\ldots i+d)$.  Then
$\theta\wedge \delta=\mu_{ij}$ and
$\theta\vee\delta=\lambda_{ij}$.
In view of Theorems \ref{sing} and \ref{no-loops-allowed}, we have

\emph{In $X_{d,n}$, $P_{\tau}$ is a smooth point if and only if for
every pair $\left(\theta,\delta\right)$ of join and meet
irreducibles, there is some $\alpha\in [\theta\wedge
\delta,\theta\vee \delta]$ such that $P_\tau(\alpha)$, the
$\alpha$-th co-ordinate of $P_\tau$ is non-zero}.

This is in fact the content of the conjecture of \cite{g-l}
(see \cite{g-l}, \S 11).
\end{remark}

This suggests for us to look at such pairs of join-meet
irreducibles in other distributive lattices and expect the
components of the singular locus of the associated Hibi toric
variety to be given by Theorem~\ref{mainth},(1) for the case of
$I_{d,n}$. This is not true in general, however, as the following
counter example shows.

\subsection{Counter Example}\label{counter-example} Let $\L$ be the interval $[ (1,3,4), (2,5,6)]$, a sublattice of $I_{3,6}$.
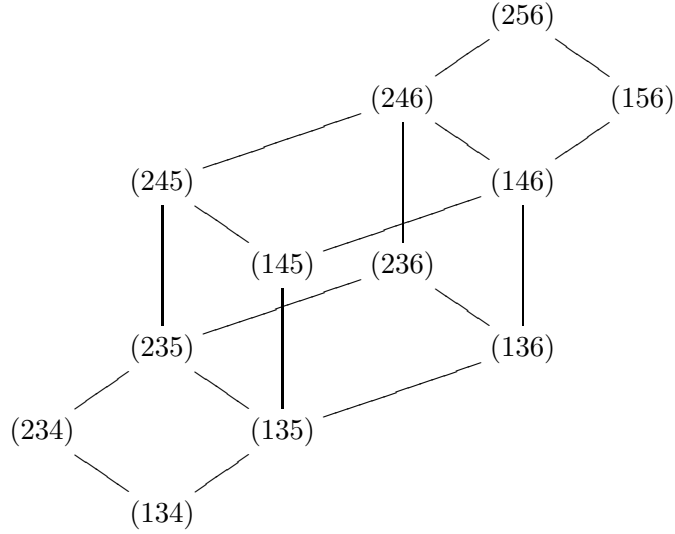
\small
\[ \xymatrix@-10pt{
 & & & & (256) & \\
 & & & (246)\ar@{-}[ur]& & (156)\ar@{-}[ul]\\
 & (245)\ar@{-}[urr]& & & (146)\ar@{-}[ur]\ar@{-}[ul] & \\
 & & (145)\ar@{-}[urr]\ar@{-}[ul]& (236)\ar@{-}[uu] & & \\
 & (235)\ar@{-}[urr]\ar@{-}[uu]& & & (136)\ar@{-}[uu]\ar@{-}[ul]& \\
(234)\ar@{-}[ur] & & (135)\ar@{-}[urr]\ar@{-}[ul]\ar@{-}[uu]& & & \\
 & (134)\ar@{-}[ur]\ar@{-}[ul] & & & &
 } \]
\normalsize

Notice that $\L$ has only one pair of join meet irreducibles,
$(2,3,4)$ and $(1,5,6)$; and thus the corresponding interval
$[\theta\wedge\delta , \theta\vee\delta ]$ is the entire lattice.
Therefore, if our result (Theorem~\ref{mainth},(1)) on the
singular locus of G-H toric varieties were to generalize to other
Hibi toric varieties, we have that any proper face is
non-singular, since any face $\tau$ must correspond to an embedded
sublattice $D_{\tau}$, and naturally this sublattice will
intersect the interval, which is just $\L$.

But this is not true!  For example, let $\tau$ be the face of
$\sigma$ such that $D_{\tau}=\{(1,5,6)\}$.  Then, the following
set
\[ \tau =C\left< e_{145}-e_{156},\,e_{136}-e_{156},\, e_{135}-e_{145},\,e_{135}-e_{136},\,e_{134}-e_{135}\right>.\]
 is a set of generators for $\tau$. Clearly $\tau$ is not generated by the subset of a basis, and thus $\tau$ is a singular face (see Lemma~\ref{basis-sing}).

Nevertheless, Theorem \ref{mainth} holds for minuscule posets as
described below.

Let $G$ be semisimple, and $P$ a maximal parabolic subgroup with
$\omega$ as the associated fundamental weight. Let $W$ (resp.
$W_P$) be the Weyl group of $G$ (resp. $P$). Then the Schubert
varieties in $G/P$ are indexed by $W/W_P$. Let $P$ be
\emph{minuscule}, i.e., the weights in the fundamental
representation associated to $\omega$ form one orbit under the
Weyl group. One knows that the Bruhat poset $W/W_P$
  of the Schubert varieties in $G/P$ is a
  distributive lattice.  See \cite{Hiller} for
  details.

  \begin{definition}\label{bruhat-def} We call
 $\mathcal{L}:=W/W_P$, a \emph{minuscule poset}, and
 $X\left(\mathcal{L}\right)$, a \emph{Bruhat-Hibi toric variety}, or a
B-H toric variety (for short).
\end{definition}

\begin{remark} Any G-H toric variety $X_{d,n}$ is a B-H toric variety.
\end{remark}

Now for $\L$ as in Definition~\ref{bruhat-def}, consider a pair
$(\alpha,\beta)$ of incomparable join-meet irreducible elements.
Let $I_{\alpha,\beta}=[(\alpha\wedge\beta ),(\alpha\vee\beta )]$.

It has been shown very recently (cf.\cite{b-l})
that a B-H toric variety $X\left(\L\right)$ is smooth at $P_\tau$
($\tau$ being a face of $\sigma$) if and only if for each
incomparable pair $(\alpha,\beta )$ of join-meet irreducibles in
$\mathcal{L}$, there exists at least one $\gamma\in
I_{\alpha,\beta}$ such that $P_\tau(\gamma)$ is non-zero.


\end{document}